\documentclass[11pt]{article}

\usepackage{amssymb, latexsym,  mathrsfs}
\usepackage{makeidx, url}
\usepackage[dvips]{graphicx,lscape} 
\usepackage{enumerate}
\usepackage[hyperref, colorlinks, linkcolor=black, citecolor=black, urlcolor=black, pagecolor=black]{hyperref}
\newif\ifpdf\ifx\pdfoutput\undefined\pdffalse\else\pdfoutput=1\pdftrue\fi
 
\newcommand{\be}{\begin{equation}}
\newcommand{\ee}{\end{equation}}
\newcommand{\ba}{\begin{array}}
\newcommand{\ea}{\end{array}}
\newcommand{\bea}{\begin{eqnarray}}
\newcommand{\eea}{\end{eqnarray}}
\newcommand{\bean}{\begin{eqnarray*}} 
\newcommand{\eean}{\end{eqnarray*}}

\newtheorem{thm}{Theorem}[section]
\newtheorem{alg}{Algorithm}[section]

\newtheorem{lem}{Lemma}[section]

\newcommand{\blem}{\begin{lem}}
\newcommand{\elem}{\end{lem}}

\newtheorem{defi}{Definition}[section]
\newcommand{\bdefi}{\begin{defi}}
\newcommand{\edefi}{\end{defi}}

\newtheorem{prop}{Proposition}[section]
\newcommand{\bprop}{\begin{prop}}
\newcommand{\eprop}{\end{prop}}

\newtheorem{coro}{Corollary}[section]
\newcommand{\bcoro}{\begin{coro}}
\newcommand{\ecoro}{\end{coro}}

\newtheorem{rmk}{Remark}[section]
\newcommand{\brmk}{\begin{rmk}}
\newcommand{\ermk}{\end{rmk}}
{
\begin{center}
 \begin{tabular}{|@{\hspace{.15in}}c@{\hspace{.15in}}|}
 \hline 
   \begin{minipage}[t]{\boxedparwidth}
   \setlength{\parindent}{.25in}}%
   {\end{minipage}  \\ 
 \hline 
 \end{tabular} 
\end{center} 
}

\newcommand{\cfpbox}[2]{\begin{center}{\fbox{\parbox{#1}{#2}}}\end{center}}

\newlength{\boxedparwidth}  

\newcommand{\stitle}[1]
{\begin{center}
   \fbox{\textbf{\large #1}}
 \end{center} }

\newcommand{\prf}{\noindent{  {\bf{\em Proof:~ }}}}

\def\sq{{{\hfill\rule{1.7mm}{1.7mm}}}}  

\newcommand{\lb}[1]{\label{#1}}
\newcommand{\rf}[1]{(\ref{#1})} 

\newcommand{\any}{\forall}
\newcommand{\orth}{\perp}    
                                
\newcommand{\bsl}{\ensuremath{ \backslash}}

\newcommand{\til}[1]{\tilde #1}

\newcommand{\lar}{\leftarrow}
\newcommand{\rar}{\rightarrow}

\newcommand{\norm}[1]{\| #1 \|}   

\newcommand{\D}[2]{\frac{\partial #1}{\partial #2}}

\newcommand{\msbm}[2]{msbm#1 scaled\magstep#2}
\newfont{\msbmtenone}{\msbm{10}{1}} 
\def\IR{\mbox {\msbmtenone R}}
\def\IC{\mbox {\msbmtenone C}}

\newcommand{\Rn}{{\IR}^n }

\newcommand{\Cn}{{\IC}^n }

\def\L1{{\cal L}_1}
\def\L2{{\cal L}_2}

\newcommand{\lam}{\lambda}

\newcommand{\eps}{\epsilon}

\newcommand{\cS}{{\cal S}}

\newcommand{\mcF}{{\mathcal F}}

\newcommand{\bA}{{\bf A}}

\newcommand{\bH}{{\bf H}}
\newcommand{\bI}{{\bf I}}

\newcommand{\bK}{{\bf K}}

\newcommand{\bM}{{\bf M}}

\newcommand{\bV}{{\bf V}}
\newcommand{\bW}{{\bf W}}

\newcommand{\bh}{{\bf h}}

\newcommand{\br}{{\bf r}}

\newcommand{\bt}{{\bf t}}

\newcommand{\bv}{{\bf v}}
\newcommand{\bw}{{\bf w}}
\newcommand{\bx}{{\bf x}}
\newcommand{\by}{{\bf y}}

\begin{document}

\title{Eigenvalue Computation from the Optimization Perspective: \\
       On Jacobi-Davidson, IIGD, RQI and Newton Updates\thanks{
This work was supported by the Mathematical, Information, and Computational Sciences Division subprogram of the Office of Advanced Scientific Computing Research, SciDAC Program, Office of Science, U.S. Department of Energy, under contract W-31-109-ENG-38.}}

\author{Yunkai Zhou\thanks{
Argonne National Laboratory, Mathematics and Computer Science Division, 
Chemistry Division, Argonne, IL 60439, U.S.A. ~Email: {\tt zhou@mcs.anl.gov}}}
\date{}
\maketitle

\begin{abstract}

We discuss the close connection between eigenvalue computation and optimization
using the Newton method and subspace methods. 
From the connection we derive a new class of Newton updates.
The new update formulation is similar to the well-known Jacobi-Davidson method.
This similarity leads to simplified versions of the Jacobi-Davidson method and
the inverse iteration generalized Davidson (IIGD) method.
We prove that the projection subspace augmented by the updating direction from
each of these methods is able to include the Rayleigh quotient iteration (RQI) 
direction. Hence, the locally quadratic (cubic for normal matrices) convergence
rate of the RQI method is retained and strengthened by the subspace methods.
The theory is supported by extensive numerical results. 
Preconditioned formulations are also briefly discussed for large scale 
eigenvalue problems.\\

{\bf Keywords:} {Eigenvalue;  Rayleigh quotient; Subspace method;  
Jacobi-Davidson; IIGD; \\Newton method} 

\end{abstract}

\section{Introduction}

We study subspace methods for the standard 
eigenvalue problem
\bea
\bA \bx & = & \lam \bx ,        \lb{seig}  
\eea 
where  $\bA$ is an  $n$ by $n$
matrix,  
$\lam$ is the 
eigenvalue and $\bx$
is the corresponding 
eigenvector.
Eigenvalue problem \rf{seig} is central to many scientific applications.
Several excellent monographs
\cite{temp00, chatel93, cul_wil02, gv96, parlet80, saad_eig92, stewar00, wilkin65}
contain in-depth discussions on properties and algorithms for eigenvalue problems,
together with many sources from science and engineering applications.

Subspace methods include the Arnoldi method \cite{arnoldi}, Lanczos method \cite{lanczos}
and the more sophisticated but important variants, 
the implicit restart Arnoldi method
\cite{sorens92} and the rational Krylov method \cite{ruhe84}. 
Other important subspace methods include the 
Davidson method \cite{davids75}, Jacobi-Davidson method \cite{jd96, gjd96, hoc_sle03},
Krylov-Schur method \cite{stewar01}, 
truncated RQ method \cite{sor_yan98, yang98}, 
and preconditioned subspace eigensolvers
\cite{cps94, fsv98, gol_ye02, knya01,kny_ney03, mor_sco86, ssf95}.

In this paper we study  mainly the Jacobi-Davidson type subspace methods for \rf{seig}. 
We highlight the connection between eigenvalue computation and optimization;
the crucial link is the Newton method for Rayleigh quotient combined with
subspace techniques. This link is well known and 
is central to the success of several eigensolvers, for example,
the Davidson method and Jacobi-Davidson method. But often 
the connection to the Newton method is emphasized, 
which is used to account for the (locally) fast
convergence property of the eigensolvers.
The importance of subspace methods is
no doubt well known. But we feel the following discussion may
still contribute to the understanding of the
Newton method applied within a subspace approach. We try to emphasize
the subspace approach.

We begin the discussion from optimization.
The traditional optimization algorithms are often
one dimensional; that is, at each iteration step, the optimization is done within 
a one dimensional space. 
As an example we look at the Newton method for 
\be \min_\bx f(\bx), ~~~  \bx \in \Rn, \lb{min}  \ee
where $f: \Rn \rar \IR$ is assumed to be twice continuously differentiable.
At step $k$  we obtain an approximate minimizer $\bx_k$. Then we solve the 
Newton correction equation:
\be  \left[ \nabla^2 f(\bx_k) \right]~ \delta \bx_k = -  \nabla f(\bx_k). \lb{newton} \ee
The approximate minimizer at step $k+1$ is updated from $\bx_k$ and $\delta \bx_k$ as
\be \bx_{k+1} := \bx_k  + \alpha_k ~ \delta \bx_k,  \lb{sigv} \ee 
where $\alpha_k$ is some line search parameter obtained by an approximate minimization;
$\alpha_k = 1$ is the standard Newton method.
Often the Armijo-Goldstein condition or Wolfe condition \cite{den_sch83} 
is used to ensure global convergence (without line search and other globalization
techniques Newton method may converge to local maximal or saddle point).
In other words, the minimizer is sought within a one dimensional space 
$span\{\delta \bx_k\}$ with affine translation $\bx_k$.

The key advantage of a subspace approach is that the convergence may be more robust and
the convergence rate may be faster than that of a single vector update approach. 
This may be seen clearly from the following: 
Instead of looking for a minimizer
within a one dimensional space $\cS_1$, we look for a minimizer 
within a multidimensional subspace  $\cS_k$ that contains   
$\cS_1$ as its subspace; hence the minimum in 
$\cS_k$ must be no greater than  the minimum in  $\cS_1$.
Again we use the Newton method for \rf{min} as an example. At step $k+1$, instead of 
looking for a single
vector update \rf{sigv}, we augment the subspace by the 
computed Newton direction (i.e.,  $span\{\bx_1, \bx_2, ..., \bx_k, \delta\bx_k\}$).
 Then
 we look for a minimizer $\til \bx_{k+1}$ within this augmented subspace. 
Then it is very clear to see 
that $f(\til \bx_{k+1}) \le f(\bx_{k+1})$. By augmenting the subspace in
this way, we 
retain the locally fast convergence property of the Newton method. 
A second advantage of the subspace approach is for clustered eigenvalues. 
Since even for the symmetric case the {\it gap theorem}
\cite{parlet80} tells us that a single vector approach (without deflation)
will not lead to a meaningful
eigenvector, we have to compute the invariant subspace related to the clustered
eigenvalues.   Another advantage of the subspace approach is that one
can seek more than one optimizer by working with different subspaces.
This feature is particularly 
important for eigenvalue computation because one often needs 
a portion of the spectrum and the related eigenvectors 
instead of a single eigenpair.


Clearly the advantage of subspace technique is not gained without price. At each
iteration of a subspace method one has to do more matrix vector products. 
And since orthonormal
basis is numerically more stable than non-orthonormal ones, one needs to  
orthonormalize the new direction against the former subspace basis.  
These imply that
the subspace dimension can not be large for
large scale problems.
To gain efficiency, one often needs to do restart and deflation \cite
{gv96, saad_eig92, sorens92, stewar00}. 
{\it Restart} is an indispensable technique  which ensures that the starting 
vector becomes progressively better, so that the convergence becomes faster 
and the subspace dimension is kept relatively small to $n$.
{\it Deflation} is also important with subspace methods. 
It tries to fix converged eigenvectors so that they will not be recomputed; and
by using only vectors orthogonal to the converged eigenvectors to augment the
projection basis, the projection
subspace will be augmented in a more efficient way. 

The important link between the eigenvalue problem \rf{seig} and optimization is
the following nonlinear function called the 
{\it Rayleigh quotient}:
\be \lb{rquo}
Q(\bx) = \frac{\bx^*\bA\bx}{\bx^*\bx}, ~~~~~~ \any   \bx  \ne {\bf 0 }.
\ee
For symmetric or Hermitian $\bA$, $\min_\bx Q(\bx)$ or  $\max_\bx Q(\bx)$ 
gives the extreme eigenvalue of $\bA$.
For general matrix $\bA$, the Rayleigh quotient defines the {\it field of value} 
\cite{greenb97, hor_joh91}:  
\[ \mcF(\bA) := \left\{Q(\bx) ~|~  \bx \in \Cn \bsl \{ {\bf 0}\} \right\}. \]
If $\bA$ is normal, then $\mcF(\bA)$ is the convex hull of the eigenvalues of  $\bA$.
In general $\mcF(\bA)$ contains the spectrum of $\bA$. With suitable subspaces one can
locate the eigenvalues of $\bA$ via Rayleigh quotient; this is the well known
{\it Rayleigh-Ritz} procedure. Certainly when the eigenvalues are located in the
complex plane, the optimization interpretation becomes obscure, since one can not
compare complex values. However, if we interpret optimizers not only 
as the end points of an interval in $\IR$ but also as the boundary of a region in $\IC$,
then  the optimization interpretation still makes sense. 
As an example of this interpretation we look at the Lanczos- or 
Arnoldi-type algorithms. It is well known that these algorithms
usually approximate the exterior eigenvalues.
Since the subspace at each iteration essentially contains the residual 
$\br_k = \bA\bx_k - \lam_k \bx_k$, for symmetric problems $\br_k$ is 
in the direction of the gradient of \rf{rquo} at $\bx_k$ (i.e., 
the steepest descent or ascent direction).  For nonsymmetric
problems $\br_k$ may not be readily considered as gradient 
since $Q(\bx)$ may not be differentiable \cite{parlet74}, 
but it is likely that the subspace augmented by $\bx_k$ contains the 
steepest descent or ascent direction (measured by vector norm) of $Q(\bx)$.
This can also be considered as the advantage of the subspace approach: one does
not have to compute the exact direction directly, but the span of the vectors contains
the best direction. In this vein of thought,
there would be no surprise that the algorithms 
converge to the exterior eigenvalues---the optimizers in the complex plane.
As for interior eigenvalues, the common approach is to apply shift and invert, that is,
to map the interior eigenvalues of the original matrix into extreme eigenvalues
of a transformed matrix, and then apply standard solvers to the transformed matrix.

{\bf Some notations:} Throughout this paper we use boldface letters
for matrices and vectors. The upper script $()^*$ denotes the transpose $()^T$
in the real case, and the conjugate transpose $()^H$ in the complex case.


\section{Subspace Eigensolvers: Davidson and Jacobi-Davidson}

The Davidson method and Jacobi-Davidson method are two of the typical examples of 
subspace methods for eigenvalue computation. Before getting into details 
we present the framework of subspace
methods for the eigenvalue problem \rf{seig} in Algorithm \ref{subs}, where
for simplicity we assume that $\bA = \bA^*$ and that only the smallest eigenvalue
and its corresponding eigenvector are sought.  Note that the optimization within
a subspace is done at step {\it 5(e)}.  

\begin{figure}[ht]
\cfpbox{143mm}{
\begin{alg} {Framework of subspace methods 
 for the eigenvalue problem \rf{seig}: }  \lb{subs} 
\begin{enumerate}
\item Start with an initial unit vector $\bx,$   $\bV_1 \lar [ \bx].$
\item Compute~ $ \bw = \bA \bx,~  \bW_1 \lar [\bw].$
\item Compute~ $ \lam = \bH_1 = \bx^* \bw, ~ \br = \bw - \lam \bx $.
\item If $\norm{\br} <= \eps$, return $(\lam, \bx)$ as the eigenpair.
\item Outer Loop:  for  $j = 1,  2, ..., m$ do
\begin{enumerate}
 \item 
       \fbox{ Call specific subspace method to construct an
        augmentation vector~ $\bt$. }
 \item Orthonormalize $\bt$ against $\bV_j$ to get a unit vector~ $\bv$.
 \item Compute~ $\bw = \bA \bv $,  ~$\bW_{j+1} \lar [\bW_j\; | \; \bw],
 ~ \bV_{j+1} \lar [\bV_j \;|\; \bv]. $
 \item Compute $ \left[ \ba{c}  \bh \\ \alpha \ea \right]  = \bV_{j+1}^*
 \bw,  ~~ \bH_{j+1}  =  \left[ \ba{cc} \bH_j & \bh \\ \bh^* & \alpha  \ea
 \right]. $
 \item Compute the smallest eigenpair $(\lam, \by)$ of $\bH_{j+1}$.
  $(\norm{\by} =1)$.
 \item Compute the Ritz vector~ $\bx = \bV_{j+1} \by$ \\
 and the residual vector ~$\br = \bA \bx -\lam \bx = \bW_{j+1} \by - \lam \bx.$
 \item Test for convergence: if  $\norm{\br} <= \eps$,  
  return $(\lam, \bx)$ as the eigenpair.
\end{enumerate}
\item Restart:  Set  $ \bV_1 = [ \bx ], ~ \bW_1 =
  [\bW_{j+1} \by], ~ \bH_1 = [\lam] $, goto step {\it 5}.
\end{enumerate}
\end{alg}

}
\end{figure}

Note also that the stopping criterion $\norm{\br} <= \eps$ is in the relative sense;
i.e., $\eps$ is related to the computed Ritz values. This is more important to
the nonnormal problems.
For $\bA \ne \bA^*$,
one can easily adapt the framework by modifying 
step {\it 5(d)} as
\[ 
 \bH_{j+1}  =   \left[ \ba{cc} \bH_j & \bV_j^* \bw \\
  \bv^* \bW_j & \bv^* \bw  \ea
 \right], \]
note only the $\bh^*$ term at step {\it 5(d)} needs to be changed 
into $ \bv^* \bW_j$.
If more eigenpairs are required, one can compute eigenpairs
 of $\bH_{j+1}$ according to 
different selection criteria   at step {\it 5(e)}.
The selection criteria include the first $l$ ($l \ge 1$) eigenvalues of the
 largest or smallest real parts, imaginary parts  or magnitudes.
Then  at step  {\it 5(a)} one computes the 
augmentation vectors from the residual vectors related to the selected
eigenpairs of $\bH_{j+1}$. 

The other important numerical step is the {\it orthogonalization} at
step {\it 5(b)}. A common choice is the modified Gram-Schmidt (MGS)
method. MGS is in practice much more stable than
Gram-Schmidt (GS), but the orthogonality of the basis constructed 
by MGS depends on the condition number of the original set of vectors.
Another problem is that MGS can not be expressed by Level-2 BLAS,
hence parallel implementation needs more communication \cite{temp00,ddsv98}.
Another  choice for orthogonalization is the DGKS correction \cite{dgks}. 
DGKS can be 
implemented more efficiently in parallel because DGKS is actually GS with 
selective reorthogonalization. To make $\bt$ orthogonal to $\bV_j$ we
do the following,
\bea 
\bt_{tmp} \lar  \bV_j^* \bt, ~~
\bv_{j+1} \lar \bt -   \bV_j \bt_{tmp}, ~~ \bv_{j+1} \lar {{\bv_{j+1}}
\over \norm{\bv_{j+1}} },  \lb{pre} \\
{\rm if}~  (\norm{\bt_{tmp}} > tol * \norm{\bv_{j+1}})~ {\rm then~ do~ the~
reorthogonalization:}  \nonumber \\
\bt_{tmp}  \lar  \bV_j^* \bv_{j+1}, ~~ \bv_{j+1} \lar \bv_{j+1} -
\bV_j \bt_{tmp}, ~~\bv_{j+1} \lar {{\bv_{j+1}}
\over \norm{\bv_{j+1}} }, 
\eea
where the criteria means that if 
the angle between $\bt_{tmp}$ and $\bv_{j+1}$ at
line \rf{pre} is smaller than $actan(tol)$, then
reorthogonalization is necessary. 

 The subspace dimension and the convergence rate are
a big concern for large scale problems. 
For subspace eigensolvers the convergence rate often depends
heavily on the quality of the starting vector. 
Restart can be used to refine the starting vector.
We know that before restart, useful information may have been 
accumulated in the recent expansion vectors. Thus it may not be desirable to keep
only the most recent
expansion vector and throw away all the rest. In \cite{ssw98, wu_sim98} 
{\it thick restart} techniques are proposed. The idea is to keep as much useful
information as possible by keeping more than one expansion vector at the first step 
of  restart.  This idea has proved to be efficient and easy to 
incorporate into the framework.

What distinguishes each subspace method is the vector augmentation at step  {\it 5(a)}.
When iterative methods (e.g., GMRES \cite{saad86}, BICGSTAB \cite{bicgstab}) 
are used to solve for the augmentation vector, step {\it 5(a)}
necessarily contains the {\it inner iteration} of the subspace methods.

The {\bf Davidson} method \cite{cps94, davids75, mor_sco86, saad_eig92} 
constructs $\bt$ as follows:
\cfpbox{137mm} {
{ ~~~ \it 5(a)}~ Construct preconditioner $\bM_j$; 
 ~ solve for~ $\bt$ from:
\[ (\bM_j - \lam \bI)\; \bt = -\br. \]
}
In the original paper by Davidson \cite{davids75}, $\bM_j = diag(\bA) $.
Originally Davidson derived this formulation by taking the derivative of $Q(\bx_j)$,
varying only the $i$-th component of $\bx_j$ (denoted as $\bx_{j(i)}$) each time:
\be  \left. \D{Q(\bx_j)}{\bx_{j(i)}} \right|_{\bx_{j(i)}+ \delta \bx_{j(i)}}  = 0,  
~~~~ i = 1, 2, ..., n.  \lb{daveq} \ee 
From \rf{daveq} we get $\delta \bx_{j(i)} = 
(\lam_j - a_{ii})^{-1} (\bA \bx_j- \lam_j \bx_j)_{(i)}$.
In \cite{davids93} Davidson related this formulation to the Newton method.
It was observed that the
convergence rate is related to how well the diagonal
matrix  $diag(\bA)$  approximates $\bA$. 
For example, if $\bA$ is diagonally dominant (in eigenvalue problems, diagonal
dominance means that the off-diagonal elements are small 
compared with the changes in magnitude
between diagonal elements \cite{wilkin65, mor_sco86, ssf95}),
 then the Davidson method converges very
fast.
Hence the diagonal matrix 
($diag(\bA) - \lam_j \bI$)
was considered as a straightforward preconditioner to $\bA -\lam_j \bI$.
The main advantage of a diagonal preconditioner is that the preconditioned
system (which is diagonal) is trivial to solve,
 but it requires $\bA$ to be  diagonally dominant.
More sophisticated preconditioners than
the diagonal matrix have been tried \cite{cps94, mor_sco86}, leading to
the so called {\it generalized Davidson method.}
As has been pointed out (e.g., in \cite{jd96}), 
the preconditioner interpretation does not
explain the improved convergence rate well since 
the exact preconditioner $(\bA -\lam_j \bI)^{-1}$ leads to stagnation 
($(\bA -\lam_j \bI)^{-1} \br$ is $\bx_j$, which lies in 
the original projection subspace).

The {\bf Jacobi-Davidson} method \cite{jd96} is an important advance for subspace
methods in eigenvalue computation. At each iteration the correction vector $\bt$
is required to reside within  the orthogonal complement of the existing projection 
subspace. Jacobi-Davidson solves the following projected equation:
\cfpbox{137mm} {
{ ~~~ \it 5(a)}~~ Solve (approximately) for ~$\bt$~ s.t. ~ $ \bt \perp \bx$ ~
 and 
\be \lb{jdeq}
 (\bI - \bx \bx^*) (\bA - \lam \bI) (\bI - \bx \bx^*) \bt = - \br.
\ee
}

%

Usually \rf{jdeq} is solved  inexactly by iterative methods.
The matrix $(\bI - \bx \bx^*) (\bA - \lam \bI) (\bI - \bx \bx^*)$ 
is always singular, but this singularity poses no essential difficulty 
to iterative methods for \rf{jdeq}.

Existing preconditioners for linear systems may be applied when solving \rf{jdeq}
iteratively. This is regarded as one of the major advantages of (Jacobi-)Davidson-type 
subspace methods, since preconditioners may not be easily incorporated into the
Arnoldi-type methods. Rational Krylov-type or (single vector) 
RQI-type methods often require rather exact system solves, 
while (Jacobi-)Davidson-type methods allow approximate system solves.

An equivalent
formulation to  \rf{jdeq} is the {\it inverse iteration generalized Davidson} (IIGD)
method by Olsen {\it et al} \cite{ojs90}.  
IIGD solves 
\be  \lb{iigd}
  (\bA - \lam \bI) \bt =  - \br  + \bx \tau,
\ee
where $\tau$ is set to ensure $\bt \perp \bx$, i.e., 
$\tau = { {\bx^*(\bA - \lam \bI)^{-1} \br } \over 
  {\bx^* (\bA - \lam \bI)^{-1} \bx} }$.  From \rf{iigd} it is clear that the expansion 
vector $\bt$ contains information in both directions $(\bA - \lam \bI)^{-1} \br$
and $(\bA - \lam \bI)^{-1} \bx$. In the original paper \cite{ojs90} the authors 
dealt with a huge problem at that time and they could not afford the subspace approach;
hence they ended up using single vector updating approach in their program 
(this is also pointed out in \cite{stewar01}). 

In \cite{jd_in} and \cite{ojs90} it was argued that Jacobi-Davidson 
and IIGD are (inexact) Newton-Raphson  
methods.  The generalized Rayleigh quotient
 $ Q(\bx,\by) := { {\bx^* \bA \by} \over {\bx^* \by} } $
was used in  \cite{jd_in} to establish the argument.
In the following section we derived a Newton update straightforwardly from the
standard Rayleigh quotient \rf{rquo}. 


\section{A New Class of Newton Update}  \lb{class}

The locally fast convergence property of the Newton method in optimization 
is an attractive feature for designing fast algorithms. Davidson,
Jacobi-Davidson and IIGD  methods all tried to establish a connection to
Newton method to explain their fast convergent behavior.  
Yet no  Newton method
combined with subspace projection applied directly to the Rayleigh quotient \rf{rquo}  
seems to have been tried.   In this section we establish this interesting application.

We assume that $\bA^* = \bA$. 
(This is not for notational simplicity
reason, since the nonnormal case may be completely different.)
The gradient of  the Rayleigh quotient \rf{rquo} is
\bea \lb{gradrq}
\nabla Q(\bx)& =& \frac{2 \bA \bx}{\bx^*\bx} - 
                  \frac{2\bx^*\bA\bx\bx}{(\bx^*\bx)^2} \nonumber\\
             &= & \frac{2}{\bx^*\bx} (\bA \bx - Q(\bx) \bx).
\eea
The Hessian of \rf{rquo} is
\bea  \lb{hessrq}
\nabla^2 Q(\bx)& =& \frac{2 \bA \bx}{\bx^*\bx} - \frac{4 \bA\bx\bx^*}{(\bx^*\bx)^2}
- \frac{2}{\bx\bx^*}\left(\bx (\nabla Q(\bx)^*) + Q(\bx) \bI\right) + 
\frac{4}{(\bx^*\bx)^2}\left(Q(\bx)\bx \bx^*\right)  \nonumber\\
  &= & \frac{2}{\bx^*\bx}( \bA - Q(\bx) \bI ) 
     - \frac{2}{\bx^*\bx} ( \bx (\nabla Q(\bx)^*) + (\nabla Q(\bx))  \bx^* ) \nonumber\\
  & =&  \frac{2}{\bx^*\bx}( \bA - Q(\bx) \bI )
      - \frac{4}{(\bx^*\bx)^2} ( \bA\bx\bx^* + \bx\bx^* \bA^* - 2Q(\bx) \bx\bx^*).
\eea
This seemingly complicated formula can be simplified. If
we assume $\bx$ is already normalized and denote $\lam = Q(x)$, then
\be \lb{hess}
\nabla^2 Q(\bx) = 2 ( \bA - \lam \bI ) - 
                  4 ( \bA\bx\bx^* + \bx\bx^* \bA^* - 2 \lam \bx\bx^*).
\ee
So the Newton equation $\nabla^2 Q(\bx) \bt = - \nabla Q(\bx)$ is
\be \lb{new} 
  [ ( \bA - \lam \bI ) - 
                  2 ( \bA\bx\bx^* + \bx\bx^* \bA^* - 2 \lam \bx\bx^*)] \bt
  = - (\bA \bx - \lam \bx) = - \br.
\ee

To see the relation with the Jacobi-Davidson equation \rf{jdeq}, we expand the
matrix in the left hand side of  \rf{jdeq} and get:  (note $ \bx^* \bx = 1$)
\bea 
&~& [ (\bI - \bx \bx^*) (\bA - \lam \bI) (\bI - \bx \bx^*)] \bt \nonumber \\
&=& [(\bA - \lam \bI) - (\bA\bx \bx^* + \bx \bx^*\bA -2 \lam \bx \bx^*)]\bt = -\br.
 \lb{jd2}
\eea


From \rf{new} and \rf{jd2} we see the essential
difference is that the scaling of the correction term 
$ (\bA\bx\bx^* + \bx\bx^* \bA^* - 2 \lam \bx\bx^*)$ to $(\bA - \lam \bI )$
in \rf{new} is twice as much as the one in \rf{jd2}. 
The formula
we derived here is Newton method, hence it preserves the locally quadratic convergence
rate. No additional treatment is needed to establish the connection to the Newton method.

Our new algorithm using Newton direction within the subspace method framework 
(Algorithm \ref{subs}) may be stated as follows:

\cfpbox{137mm} {
{ ~~~ \it 5(a)}~~ Solve (approximately) the following equation for $\bt \orth \bx$:
\be \lb{nt}
  (~ \bA - \lam \bI 
      -2( \bA\bx\bx^* + \bx\bx^* \bA^* - 2 \lam \bx\bx^*)~ )  \bt = - \br 
\ee
}

In nonlinear optimization, the Hessian matrix is often difficult to obtain,
hence low rank updates such as DFP and BFGS are used to obtain approximate Hessian
matrix, resulting in the Quasi-Newton methods with a super linear convergence rate.
 For eigenvalue computation, the situation is better because we can compute 
the explicit Hessian of the Rayleigh quotient. Actually in \cite{wss98} the correct
Hessian was derived, but the formula in \cite{wss98} was not simplified and
it was more complicated than \rf{hessrq};
hence the authors did not try applying Newton method to the Rayleigh Quotient, 
while they proposed other very neat Newton-type methods based on transformed 
formulations of the eigenvalue problem.

The other difference between  \rf{new} and the Jacobi-Davidson equation \rf{jd2} 
is that the $\bx\bx^* \bA^*$ term in \rf{nt} is 
$\bx\bx^* \bA$ in \rf{jd2}. This certainly makes no difference under 
the assumption $\bA=\bA^*$. Actually under this assumption, equation \rf{nt}
is equivalent to:
\be  \lb{new2}
 (\bI - 2 \bx \bx^*) (\bA - \lam \bI) (\bI - 2 \bx \bx^*) \bt = - \br.
\ee
This looks very similar to the Jacobi-Davidson equation \rf{jdeq}. We 
note that \rf{jd2} was derived from computational insight (where
the idea of the Jacobi orthogonal component correction is exploited as the
orthogonal projector). We will discuss later that this is the truly valuable insight
that makes the Jacobi-Davidson to be so efficient.
While \rf{new2} is derived
rigorously by the Newton method. The major difference is that  
$(\bI - 2 \bx \bx^*) $ is perfectly conditioned since $ (\bI - 2 \bx \bx^*)^2 = \bI$
($\bI - 2 \bx \bx^* $ is also called the {\it Householder reflector}),
while the orthogonal projector $(\bI - \bx \bx^*)$ is singular.
Hence, away from an eigenvalue, 
\rf{new2} is better conditioned than \rf{jd2}.
When near an eigenvalue $\lam$, both \rf{jd2} and \rf{new2} relate to
the inverse iteration on $(\bA - \lam \bI)$.  
As discussed in \cite{parlet80, pet_wil79}, the ill-conditioning of $(\bA - \lam \bI)$ 
does not
prevent computing the direction of the eigenvector for $\lam$ correctly.

For large scale problems, computing the Newton direction exactly is usually not practical.
Equation \rf{new2} need to be replaced by some preconditioned equation.
We can incorporate preconditioners in the same way as
the generalized Davidson method or Jacobi-Davidson method does.
 I.e., if we construct a 
good preconditioner $\bM_j$  at step $j$, then we solve the following:
\cfpbox{137mm} {
{ ~~~ \it 5(a)}~~ Applying preconditioner $\bM_j$, 
 solve the following equation for $\bt \orth \bx$:
\be  \lb{new3}
 (\bI - 2 \bx \bx^*) (\bM_j - \lam \bI) (\bI - 2 \bx \bx^*) \bt = - \br.
\ee
}
If the formulation in \rf{nt} is preferred, then we solve for $\bt \orth \bx$
from 
\be \lb{newton2}
 (~ \bM_j - \lam \bI 
      -2( \bA\bx\bx^* + \bx\bx^* \bA^* - 2 \lam \bx\bx^*)~ )  \bt = - \br.
\ee

In both formulations \rf{new3} and \rf{newton2}, 
we may apply the easily available diagonal preconditioner 
$\bM_j=diag(\bA)$.  We note that they now can not be solved as trivially
as diagonal systems.  We also note that \rf{newton2} contains more information 
about the original $\bA$ than does \rf{new3}.

The overhead caused by the correction terms in \rf{new3} or \rf{newton2} is marginal
for iterative methods because
 one mainly computes the approximate solution by matrix-vector
products.
But we point out that the condition number of the matrix (which need not be formed 
explicitly) affects the convergence rate of any iterative method.
For eigensolvers that utilize preconditioners from linear systems directly,
another complication may occur, as pointed out in \cite{sorens02}, that a good 
preconditioner $\bM_j$ for $\bA$ (such as preconditioners from multigrid method) 
need not necessary be a good preconditioner for
$\bA - \lam \bI$ or other modifications to $\bA$. 
We also note that the interpretation of conditioning is different
between linear systems and eigenvalue problems: for linear system, the conditioning
 depends solely on $\sigma_{max}/\sigma_{min}$ where $\sigma$ denotes the singular 
value; while for eigenvalue problems, 
each eigenvalue has its own conditioner number and the separation of eigenvalues
is closely related to the convergence rate of eigensolvers.

Note that $\br = \bA \bx - \lam \bx$ is
computed at each step of the Davidson-type methods. 
If $\lam \in \IR$, the three correction terms
in \rf{nt} and \rf{newton2} may be simplified into two terms:
\be \lb{simp}
\bA\bx\bx^* + \bx\bx^* \bA^* - 2 \lam \bx\bx^* = \br \bx^* + \bx\br^*.
\ee
This simplification makes updating the matrix vector products at each iteration easier,
and we use it in our code for the numerical tests.


\section{Unification of Jacobi-Davidson, IIGD and Newton Updates via RQI} \lb{unify}

In large scale eigenvalue computation, solving preconditioned system (inexactly)
via iterative methods is a crucial part for most of the fast eigensolvers.
We see from Section \ref{class} that when suitable correction equations are found,
preconditioners may be
readily adapted into the existing subspace methods by modifying the correction
equations. Hence in the rest of this paper we
concentrate on the derivation of the correction equations, without 
being distracted by preconditioners.

When we talk about fast convergent algorithms for eigenvalue problems,
an important candidate that comes to mind is the Rayleigh quotient iteration 
(RQI). RQI is 
cubic convergent for normal matrices and 
quadratic
convergent for nonnormal matrices \cite{parlet74}.
Extensive studies of RQ and RQI exist.  
Properties of RQI and 
different type of generalizations may be found, for example,
 in \cite{dax03, notay03, ostrow, parlet74, parlet80, stewar00} 
and references therein.

In the following we unify the Jacobi-Davidson, IIGD and the proposed methods in 
Section \ref{class} by RQI. This unification leads us to propose a simplified version
of the Jacobi-Davidson method and a simplified IIGD method.
We first establish the following theorem.
\begin{thm} \lb{rqi_nw}
Let the initial unit vector be the same $\bx$, assuming the Rayleigh quotient $\lam
= \bx^*\bA\bx$ is not an eigenvalue of $\bA$ yet. Then
the subspace method with equation \rf{jdeq} or \rf{iigd}
or \rf{new2} each produces a subspace that 
includes the RQI direction during the next iteration.
\end{thm}

\prf 
The subspace method with updating equation  \rf{jdeq} or \rf{iigd} or \rf{new2} 
is a special case of the following (with suitable $\alpha$ and $\beta$):
\cfpbox{137mm} {
{ ~~~ \it 5(a)}~
 Solve (approximately) the following equation for $\bt \orth \bx$:
\be \lb{tmp2}
(\bI -  \alpha\; \bx \bx^*) (\bA - \lam \bI) (\bI -  \beta \;\bx \bx^*) \bt = - \br,
 ~~~~~  \any~ \alpha  \ne 0.
\ee
}
We  now prove the following stronger statement:
Solving \rf{tmp2} for $\bt \orth \bx$, then the subspace augmented by $\bt$ will
include the RQI direction during the next iteration.

It is well known that during the next iteration, RQI produces the direction
\[ \bx_R = (\bA - \lam \bI)^{-1} \bx.  \]

Under the restriction  $\bx^*\bt =0$, equation \rf{tmp2} becomes
\be \lb{tmp3}
    (\bA - \lam \bI) \bt = -\br + \eta \;\bx. 
\ee
This is actually the IIGD equation.  From \rf{tmp3}
\be  \lb{update}
\bt =  (\bA - \lam \bI)^{-1}\br +  \eta\; (\bA - \lam \bI)^{-1} \bx 
 =  - \bx + \eta\; \bx_R,
\ee  where $\eta = 1 / (\bx^*\bx_R)$ is chosen 
to ensure  $\bx^*\bt =0$. 
Note that $\bx$ is in the original projection subspace, say, $\bV_j$. We see that 
the RQI direction $\bx_R$ will be included in the augmented projection subspace
$\bV_{j+1}$ when we carry out {\it 5(b)} in Algorithm \ref{subs}. 
\sq

This theorem shows that the subspace method with \rf{tmp2} is able to retain
the desirable fast convergence property of RQI.  As discussed in the introduction
section, the subspace approach possibly converges faster than RQI,
and it also possesses the other advantages of a subspace method over a 
single-vector method.

Formula \rf{tmp2} may be linked to the Newton method in other ways.
From \rf{tmp2} and  $\bx^*\bt =0$  we see
\[ 
(\bI -  \alpha\; \bx \bx^*) (\bA - \lam \bI) \bt = - \br.  
\]
This leads to:
\be  \lb{tmp4} 
(\bA - \lam \bI) \bt -  \eta\; \bx = - \br,   ~~~~ \bx^*\bt =0. 
\ee
Equation \rf{tmp4} is equivalent to the following bordered equation:
\be  \lb{tmp5} 
\left[ \ba{cc} \bA - \lam \bI  &   -\bx \\
                -\bx^*          &   0 \ea
\right] 
\left[ \ba{c} \bt \\  \eta \ea \right]  = 
\left[ \ba{c} -\br \\  0 \ea \right].
\ee
This bordered equation also appears in \cite{pet_wil79, wss98} where Newton
method is applied to the following nonlinear equation:
\be \lb{tmp6}
\left\{ \ba{r}
(\bA - \lam \bI) \bx  =  {\bf 0}, \\
 -{1\over 2}\bx^*\bx + {1\over 2} = 0.
\ea
\right.
\ee
The bordered matrix turns out to be the Jacobian matrix of \rf{tmp6}. The
Newton correction equation to \rf{tmp6} is exactly \rf{tmp5}. From this
point of view, each of the algorithms --- Jacobi-Davidson, IIGD and the generalized
form \rf{tmp2} --- can be viewed as a Newton method for \rf{tmp6} 
carried out with subspace acceleration. 
A hidden reason for the success of the  Jacobi-Davidson and IIGD methods
may be that the bordered matrix $\left[ \ba{cc} \bA - \lam \bI  &   -\bx \\
                -\bx^*          &   0 \ea
\right]$ is generally nonsingular (even when $\lam$ is a simple eigenvalue of $\bA$).
The nonsingularity was first proved in  \cite{pet_wil79}.  
In actual computations one does not need
to form the bordered matrix.  Discussion of the preconditioned
form of \rf{tmp5} may be found in \cite{sorens02}.

We prefer unifying the methods by RQI instead of by the Newton method,
since the convergence rate for RQI is locally cubic for normal matrices and 
locally quadratic for nonnormal matrices. 
Whereas the Newton method usually explains only the locally quadratic convergence 
rate for both cases.

As seen from the proof, the fast convergence rate of these
methods depends  mainly on including the RQI direction to augment the
projection subspace. This fact is well known for the Jacobi-Davidson method.
But to cast it in a more general setting as we did here is new.  
Also, the proof leads to other interesting observations, which we discuss
in the next section.


\section{Simplification of Jacobi-Davidson and IIGD}  \lb{simple}

A more critical look at the proof of Theorem \ref{rqi_nw} reveals
that the constants before the correction terms in \rf{jdeq} and \rf{new2}
(hence also  \rf{nt})
 surprisingly seem not essential in establishing the link to RQI. 
As seen from \rf{tmp2}, any nonzero $\alpha$
seems to work.  The arbitrary choice of the parameters ($\alpha, \beta$)
in \rf{tmp2} may be
convenient for algorithm design,
but certainly it is not mathematically pleasant.

The difference may be that from \rf{tmp2} to \rf{tmp3} we require
both $\bx^*\bt =0$ and implicitly 
$\eta = \alpha \;\bx^*(\bA - \lam \bI)\bt =  1 / (\bx^*\bx_R) $.
This implies that different choices of $\alpha$ (and possibly $\beta$) may not
be equivalent. But this explanation does not go further.

We have determined, through extensive numerical experiments, that $\alpha = 1$ 
in \rf{tmp2} is the most essential character that leads to the success of
the Jacobi-Davidson method.
In the Jacobi-Davidson method \cite{jd96}, $\alpha = \beta =1$ is inspired by
the Jacobi orthogonal component correction. This construction 
is indeed of deep computational insight.
While in \rf{new2}, $\alpha = \beta = 2$ is derived from the Newton method for the
Rayleigh quotient of $\bA(=\bA^*)$. 
Even though $(\bI -  2 \bx \bx^*)$ is a perfectly conditioned matrix,
without the additional restriction $\bx^*\bt =0$,  \rf{nt} and \rf{new2} 
generally converge more slowly than do the Jacobi-Davidson and IIGD.
At first sight this may seem questionable because from Theorem \ref{rqi_nw}
all the methods should be equivalent. The explanation is that multiplying
$(\bA - \lam \bI) \bt = -\br$ by perfectly conditioned matrices may not
move the solution $\bt$ away from the direction of $\bx$; hence the stagnation
for the perfectly preconditioned Davidson method may also happen to \rf{nt}
and \rf{new2}, though not as seriously. 
Additionally requiring $\bt \orth \bx$, if done after solving the correction
equations \rf{nt} or \rf{new2} by direct methods,  may lead to a poor quality $\bt$.
In this case, often the DGKS method at step {\it 5(b)} 
is not enough to ensure orthogonality. One would have to perform three steps of
Gram-Schmidt orthogonalization, which often results in a vector that 
does not augment the projection subspace efficiently.
If the  $\bt \orth \bx$ is enforced at each step of any iterative methods for  
\rf{nt} or \rf{new2}, then from the equality $(\bI - \bx\bx^*)(\bI - \gamma \bx\bx^*)
= (\bI - \gamma  \bx\bx^*)(\bI - \bx\bx^*) = (\bI - \bx\bx^*)$, for $\any \gamma$, 
we see that 
\rf{nt} and \rf{new2} are essentially Jacobi-Davidson.


As seen from Algorithm \ref{subs}, after solving
the correction equation for $\bt$ in step {\it 5(a)}, this $\bt$ is orthonormalized
against a subspace that contains $\bx$. Hence, theoretically speaking, requiring
$\bt \orth \bx$ at step {\it 5(a)} is not necessary; 
but for numerical reasons, especially when 
iterative methods are used to solve the correction equation, it may be necessary
to orthogonalize the solution $\bt$ against $\bx$ once in a while or at every step
of the inner iteration (as suggested in \cite{jd96}). 

We found out that choosing $\alpha = 1$ in \rf{tmp2} is crucial in moving
$\bt$ sufficiently away from the direction of $\bx$. 
While the value of $\beta$ does not appear to be important, it may be set to 
any value not too large as to lead to overflow. For 
cost-effective reasons, we may set $\beta = 0$ and simplify the Jacobi-Davidson method 
as follows:
\cfpbox{137mm} {
{ ~~~ \it 5(a)}~~ Solve (approximately) for ~$\bt$~ s.t. ~ $ \bt \perp \bx$ ~
 and: 
\be \lb{jdmeq}
 (\bI - \bx \bx^*) (\bA - \lam \bI)\; \bt = - \br.
\ee
}
Note that \rf{jdmeq} implies
\[   (\bI - \bx \bx^*) (\bA - \lam \bI)\; \bt = - (\bI - \bx \bx^*) \br, \]
so we actually solve the equation from a right hand side vector orthogonal to $\bx$.

For $\bA \ne \bA^*$, the unsymmetric matrix in \rf{jdmeq} may not cause 
objection. 
For  $\bA = \bA^*$, it may suggest that the symmetric matrix in \rf{jdeq}
would be preferred. 
We show by extensive numerical experiments that
for both the symmetric and unsymmetric cases, if direct methods are used to solve
\rf{jdmeq} and \rf{jdeq}, then the  cheaper
\rf{jdmeq} performs as well as \rf{jdeq}.
We point out that if iterative methods are used to solve \rf{jdmeq} and the
solution at each step is orthogonalized to $\bx$, then \rf{jdmeq} is the same 
as the original Jacobi-Davidson method. 
The excellent numerical behavior of direct methods for \rf{jdmeq}
suggests that in iterative methods, the orthogonalization of the solution to
$\bx$ may be carried out after several iteration steps instead of after each step.

One advantage of Jacobi-Davidson over IIGD is that IIGD requires two system-solves
per inner iteration, 
which is more expensive than the one system-solve Jacobi-Davidson.

Also from \rf{update} in the proof of Theorem \ref{rqi_nw}, we see that in order to
capture the fast convergence property of RQI  and to 
move $\bt$ away from the $\bx$ direction, the $\eta$ should not be too small.
And \rf{update} suggests the following simplified IIGD, which requires only one
system-solve:
\cfpbox{137mm} {
{ ~~~ \it 5(a)}~~ Solve (approximately) for $\bt$ from:
\be \lb{rqieq}
 (\bA - \lam \bI)\; \bt = \bx,
\ee
~~~ Compute the constant $\eta =  1 / (\bx^*\bt)$,~  then  assign 
$\bt \lar  (\eta \bt - \bx).$
}

This modification is simple. Other researchers may have noticed
similar equivalent formulas several years ago, but to our knowledge,
to explicitly propose solving the IIGD equation by this one-solve approach 
is new, and no numerical results seem have been published that compare the
one-solve approach with the mainly used two-solve approach.
Our modification actually saves almost half the computational
cost at each inner iteration of IIGD. Even with one less system-solve, 
it produces a $\bt$ that is readily orthogonal to $\bx$ (even when 
$(\bA - \lam \bI)^{-1} \bx$ makes a small angle with $\bx$.
This is an improvement on the observation 4.1.(c) in \cite{jd96}).
Both the above analysis and the extensive numerical
experiments that we performed show that \rf{rqieq} works as well as the IIGD.  
Because of the reduced 
(often near-singular) system solves,  \rf{rqieq} can be expected to perform better
than \rf{iigd}. 
Most important, it is the least expensive method but equally efficient to 
Jacobi-Davidson and IIGD.  

We note that existing preconditioned methods for linear systems may be easily
incorporated into \rf{jdmeq} and \rf{rqieq}.  We use \rf{rqieq} as an example:
If at step $j$ we have a good
preconditioner $\bK_j$ for $(\bA - \lam_j \bI)$, then we solve $\bK_j  \bt = \bx_j$ and
compute $\eta = { 1 \over {x_j^*t} }$.  Then the updating vector can be obtained by
$ (\eta \bt - \bx_j),$ which is always orthogonal to $\bx_j$  --- the direction we
are interested in to augment the projection basis. The significance of
\rf{rqieq} may be seen from the following: If we stick to the preconditioned form
of \rf{iigd}~ (i.e.,  $ \bK_j \bt = - \br_j + \eta \bx_j$), then 
two preconditioned equations would have to be solved 
for the orthogonalization constant $\eta$ ($\eta = {{ \bx_j^* \bK_j^{-1} \br_j } \over
 {\bx_j^* \bK_j^{-1} \bx_j }}$) since $\bK_j^{-1} \br_j$ may not be equal to $\bx_j$. 
But
our modified scheme \rf{rqieq} aptly reduced the two system-solves into one 
system-solve, even for the preconditioned formulation.
  

\section{Numerical Results and Discussions}

The purpose of this paper is not to present actual implementations of
solvers for large scale eigenvalue problems. Our purpose is to study the efficiency
of different correction equations. Different preconditioning
techniques may be applied readily after the correction equations are derived.
Currently one of the most efficient and well known eigensolvers is the 
Jacobi-Davidson method \rf{jdeq}.  IIGD \rf{iigd} is equivalent to Jacobi-Davidson and
is well known in the computational chemistry community. We propose 
the generalized form \rf{tmp2}. In this section we present numerical behavior of
the more specific forms \rf{nt}, \rf{new2}, \rf{jdmeq} and \rf{rqieq}, with comparisons
among themselves and to the well known \rf{jdeq} and \rf{iigd}.
In the result-plot Figures, 
N1 and N2  denote \rf{nt} and \rf{new2}, respectively; JD and IIGD 
denote \rf{jdeq} and \rf{iigd}, respectively; 
and JDm and IIGDm denote the modified forms
\rf{jdmeq} and \rf{rqieq}, respectively.

The numerical experiments were done by using Matlab v6.5 on Pentium IV PCs with OS
Linux Redhat7.3.  The models are restricted to relatively small scale
problems ($200 \le n \le 1000$). 
We use direct methods (Matlab ``\bsl") to solve the correction equations,
since for matrices of this scale 
``\bsl"  is often at least ten times faster and more accurate
than iterative methods.  We note that when $\bM$ is singular,  $\bM \bsl \br$
give the solution $\bM^{+}\br$ where  $\bM^{+}$ is the pseudo-inverse of  $\bM$.

As pointed out in Section \ref{simple}, in Algorithm \ref{subs},
the solution $\bt$ from step {\it 5(a)} 
is made orthogonal to $\bx$ in step {\it 5(b)}. Hence, for the Jacobi-Davidson
type equations \rf{jdeq}, \rf{nt}, \rf{new2} and \rf{jdmeq}, 
we {\bf did not} additionally require $\bt \orth \bx$ in step  {\it 5(a)}. This 
modification will 
tell which equation would {\it naturally} provide a solution that can 
augment the subspace efficiently.  

We tried the methods on models from the Matrix Market.\footnote{
\url{http://math.nist.gov/MatrixMarket/}}
Models from this website are often related to 
realistic and 
difficult problems.
For all the tests we chose relative tolerance $tol = 10^{-10}$ as the
stopping measure for the residual norm.
We report results for the models: {\tt bwm200, ck400, pde225, pde900, rdb200, rdb450,
rdb968, rw496}.  The name of each model is listed in the first term at the title of each
subgraph. These models may be retrieved easily by searching for 
eigenvalue problems in Matlab format at the Matrix Market website.
The physical meaning
of these models may be found at the same website. The 
{\tt mode} denotes the eigenvalue selection criteria, 
where LR, LM, SR and SM stand for largest real part,
largest magnitude, smallest real part and smallest magnitude, respectively. The modes
were chosen according to the description of each model. 
Since there are fewer symmetric standard eigenvalue problems in Matlab
format at Matrix Market, we also applied the methods to matrices $(\bA+\bA^*)/2$ 
when $\bA \ne \bA^*$. We note that this simple symmetrization still led to problems 
that are more realistic (and usually more difficult) than randomly 
constructed models.

The methods discussed
in this paper generally converge fast locally, but not globally.  
For the problems from Matrix Market, without a good initial vector, 
each method may take a long time to arrive at the fast convergent region. 
Hence, we chose the initial vector as follows.
We computed the eigenvector corresponding to the required eigenvalue
by the Matlab function ${eig}$ and perturbed it as the initial vector. 
For the nonsymmetric problems, we perturbed
the eigenvector by  $10^{-2}*rand(n,1)$; for symmetric problems, the perturbation was
set as $5*10^{-2}*rand(n,1)$.  By choosing the initial vector this way, the 
locally fast convergence behavior of each method becomes clear to observe and easy to
compare. Certainly in real applications it is not easy to find a good initial guess
without a good understanding of the models at hand,
and one may have to resort to restart techniques. The approach we took in this section
was mainly to verify the analysis we did in this paper, as well as to 
check the viability of the proposed formulations. 

As seen from the numerical results,
the Jacobi-Davidson, IIDG and their simplified versions
\rf{jdmeq} \rf{rqieq} do show locally quadratic convergence for nonsymmetric problems
and cubic convergence for symmetric problems. And the convergence appeared faster 
than the theoretic rate in some
cases, which is due to the subspace acceleration. For all the problems tried, 
JDm and IIGDm  behaved almost the same as Jacobi-Davidson and IIGD; 
as seen from the plots, the four lines
are often indistinguishable.  This result also agrees with our analysis
in Sections \ref{unify} and \ref{simple}.

What needs to be explained more is the behavior of \rf{nt} and \rf{new2}. 
Although for some
constructed models \rf{nt} and \rf{new2} appeared to be as competitive as JD and IIGD,
for most of the more realistic models they performed more slowly. The observed 
convergence behavior for \rf{nt} and \rf{new2} is locally quadratic for symmetric
problems and linear for nonsymmetric problems. It may be super-quadratic 
or super-linear for each case because of the subspace acceleration. The reason is that
the derivation of the gradient and Hessian in Section \ref{class} holds true for 
$\bA=\bA^*$, whereas
the gradient for the Rayleigh quotient of $\bA \ne \bA^*$ is different and may not exist. 
We point out that even though in the Matlab codes we did not require 
$\bt \orth \bx$ at step {\it 5(a)} for \rf{nt} and \rf{new2},
adding this requirement
may not help much if the correction equations are solved by direct methods.
The reason is discussed in Section \ref{simple}. 
If this requirement is enforced at each step of an iterative method
for  \rf{nt} and \rf{new2},  then  \rf{nt} and \rf{new2} essentially become 
Jacobi-Davidson and IIGD, as seen from the proof of Theorem \ref{rqi_nw}.

The
Jacobi-Davidson equation \rf{jdeq} together with equations \rf{nt} and \rf{new2} 
led to the generalized formulation \rf{tmp2},
which led to the unified theory
and the two simplifications \rf{jdmeq} and \rf{rqieq}.  The modifications are
cheaper and simpler than the original formulations. Analysis and numerical evidences 
show that they are as efficient as the original formulations.  We expect that
preconditioned formulations based on \rf{jdmeq} and \rf{rqieq} will be practical in
actual applications.

The numerical results we present here evidently show that 
$\alpha = 1$ in \rf{tmp2}  is fundamental 
for the efficiency of the methods we discussed in this paper.
We hope that they also contribute to the understanding and appreciation of the 
Jacobi-Davidson method and the IIGD method.

One problem related to Jacobi-Davidson and IIGD type methods
is the possible slow convergence when a starting vector is poor. Our unreported 
numerical tests show that with a random initial vector 
many iterations may be needed to reach the fast convergent region
(some of the numerical results in \cite{hoc_sle03} showed similar behavior).
This suggests that efficient globalization techniques may need to be developed. 
A practical approach is to apply a few steps of cheaper Lanczos or Arnoldi 
iteration to get a better starting vector and then integrate restart during the
Davidson type outer iteration.  Another approach is to modify shifts
in the correction equations.
This may require a better understanding on the difference between 
preconditioning for linear systems and preconditioning for eigenvalue problems.
\cite{ssf95} contains very interesting results in this direction
for the symmetric eigenvalue problem. 

%


\begin{figure}[htp]
\begin{center}
\scalebox{.397}{\includegraphics{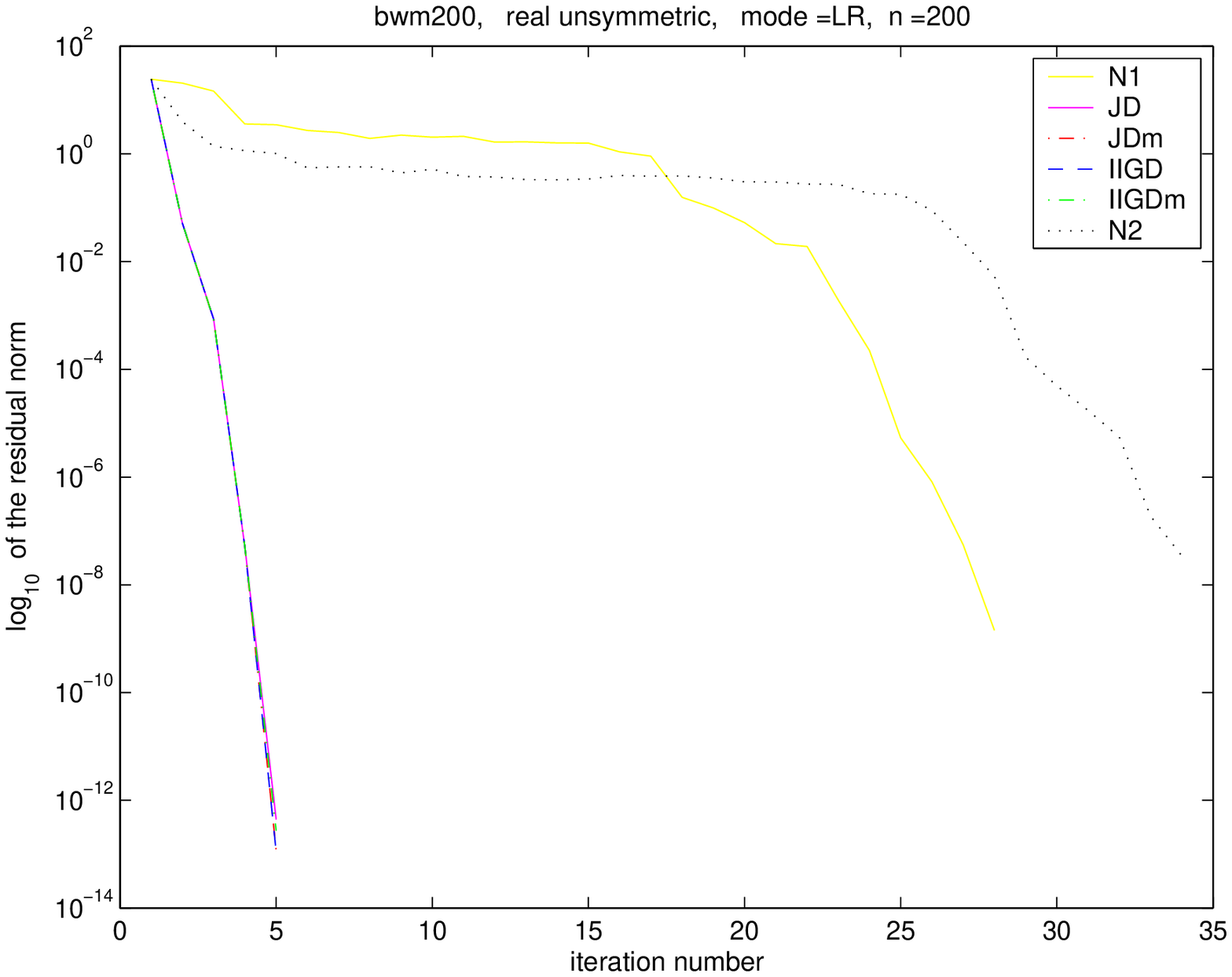}}
\scalebox{.397}{\includegraphics{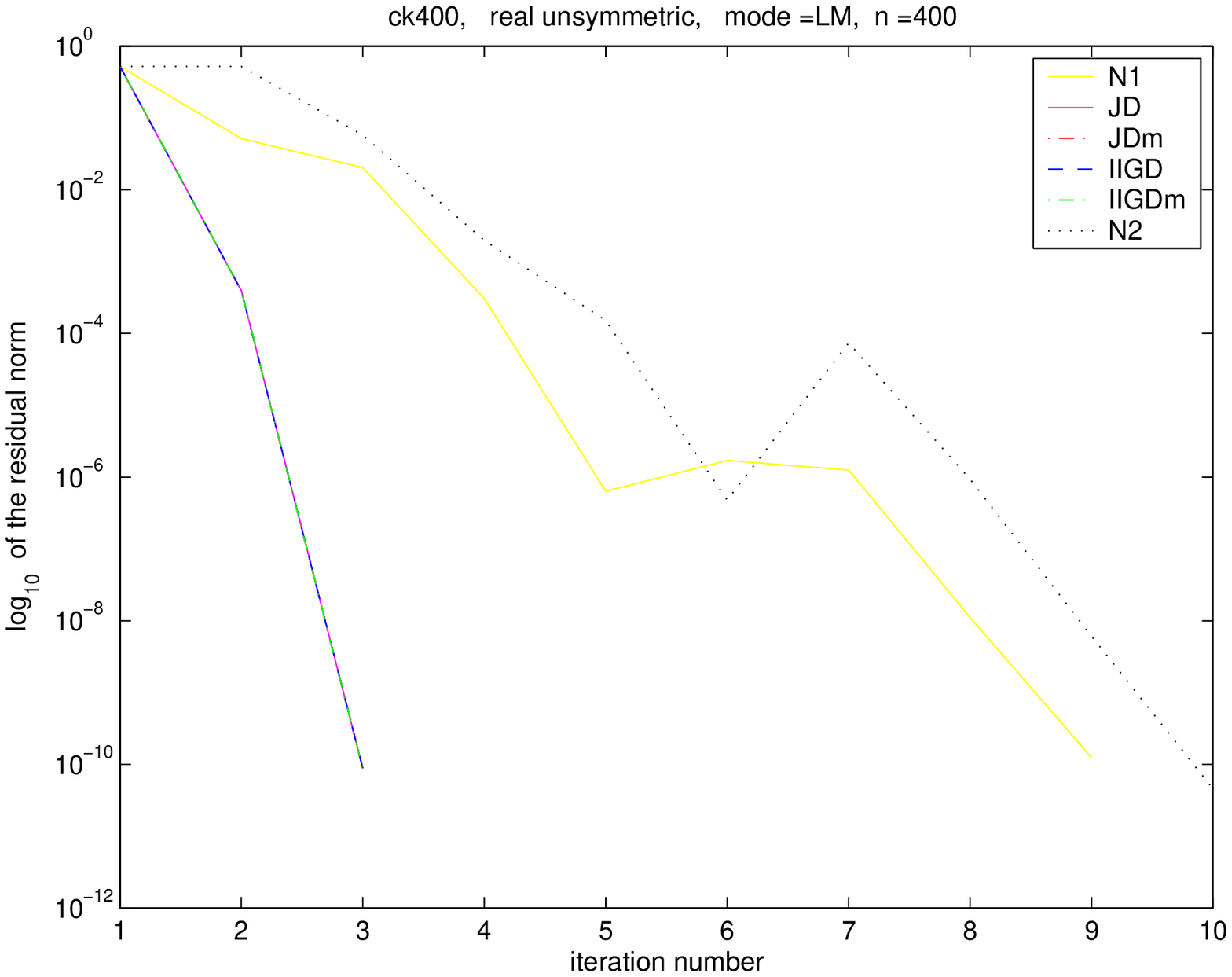}}

\vspace*{.2cm}

\scalebox{.397}{\includegraphics{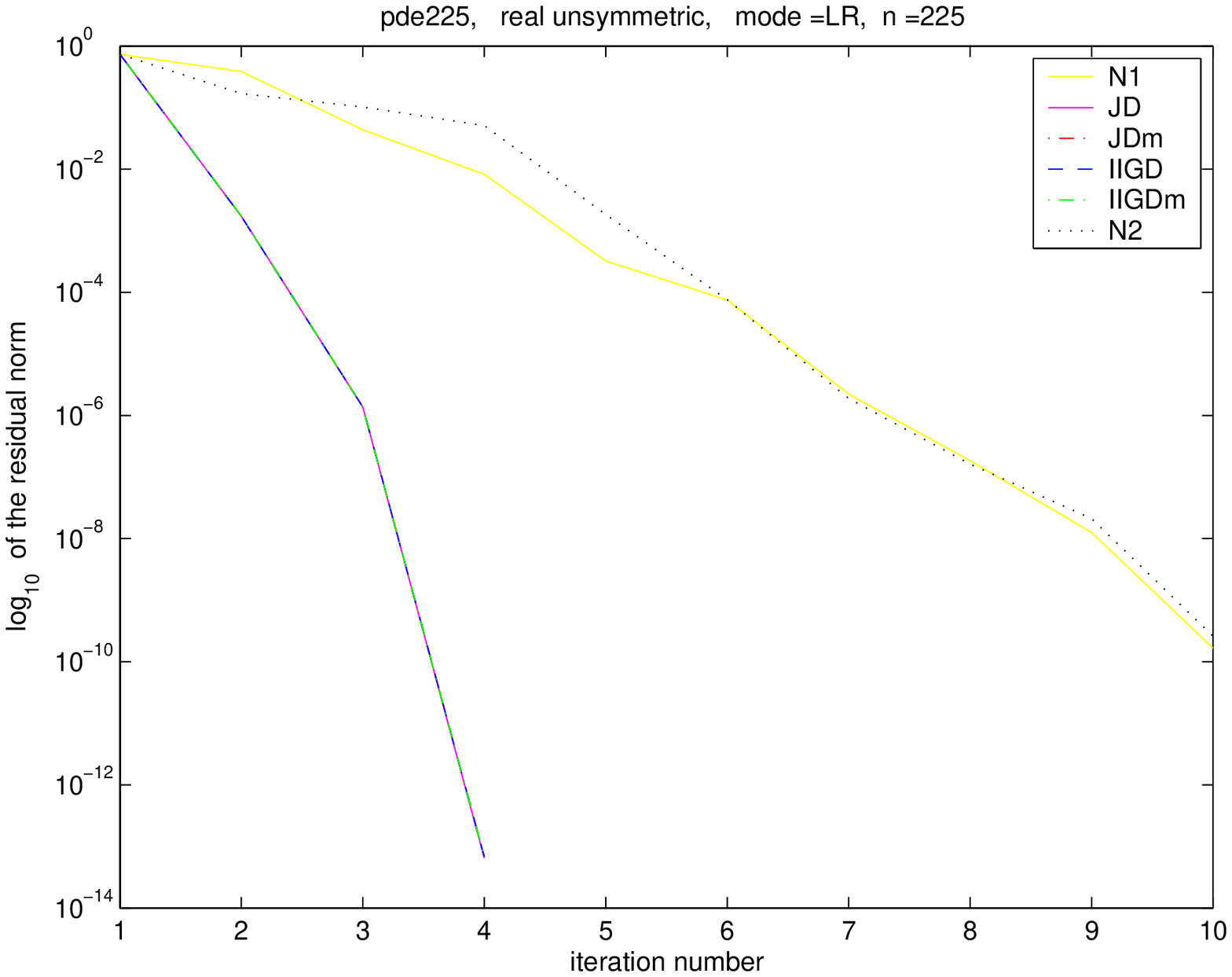}}
\scalebox{.397}{\includegraphics{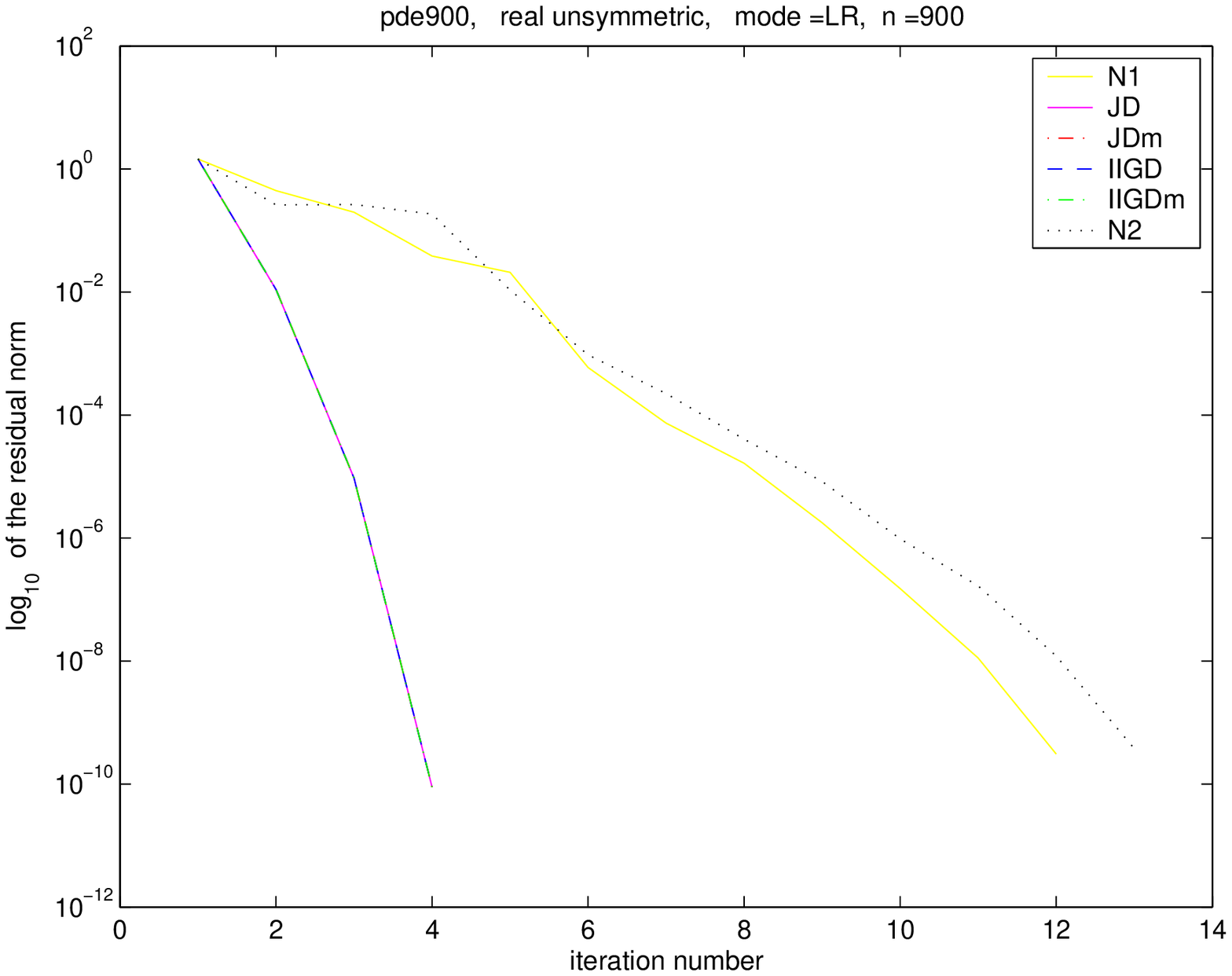}}

\end{center}
\caption{Comparison of performance for nonsymmetric models: 
bwm200, ck400, pde225, pde900. }
\label{nonsym1}
\end{figure}

\begin{figure}[ht]
\begin{center}

\scalebox{.397}{\includegraphics{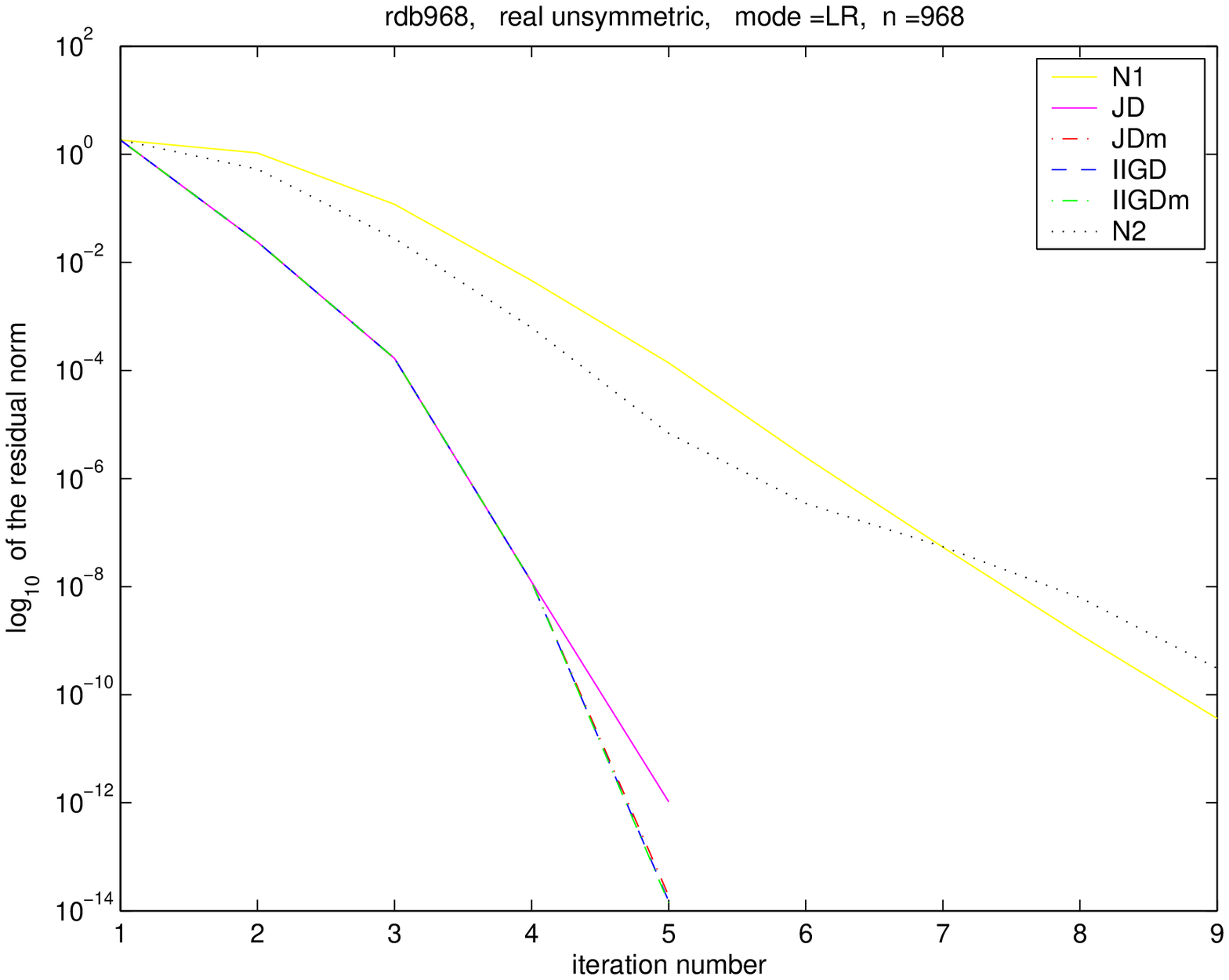}}
\scalebox{.397}{\includegraphics{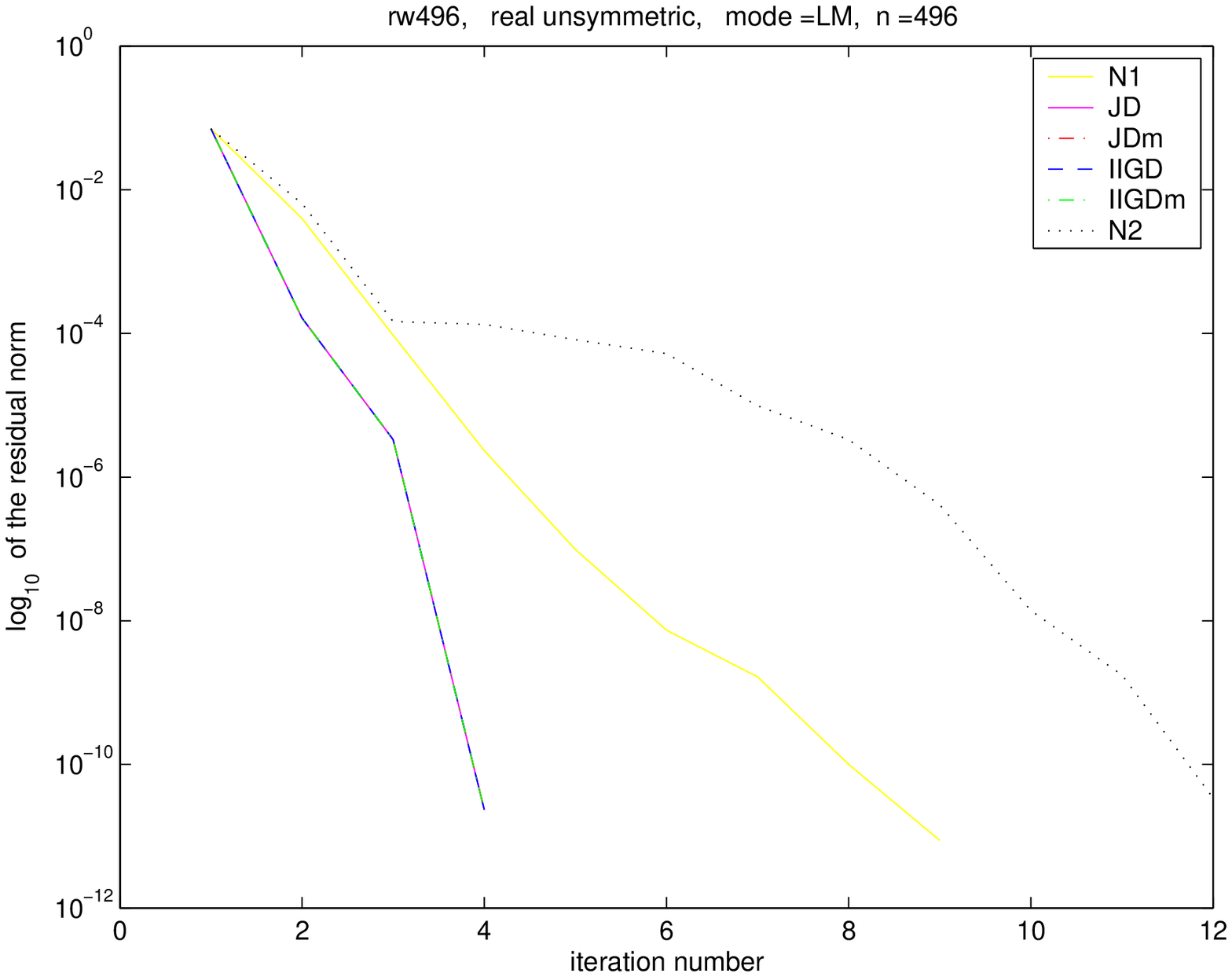}}

\vspace*{.2cm}
\scalebox{.397}{\includegraphics{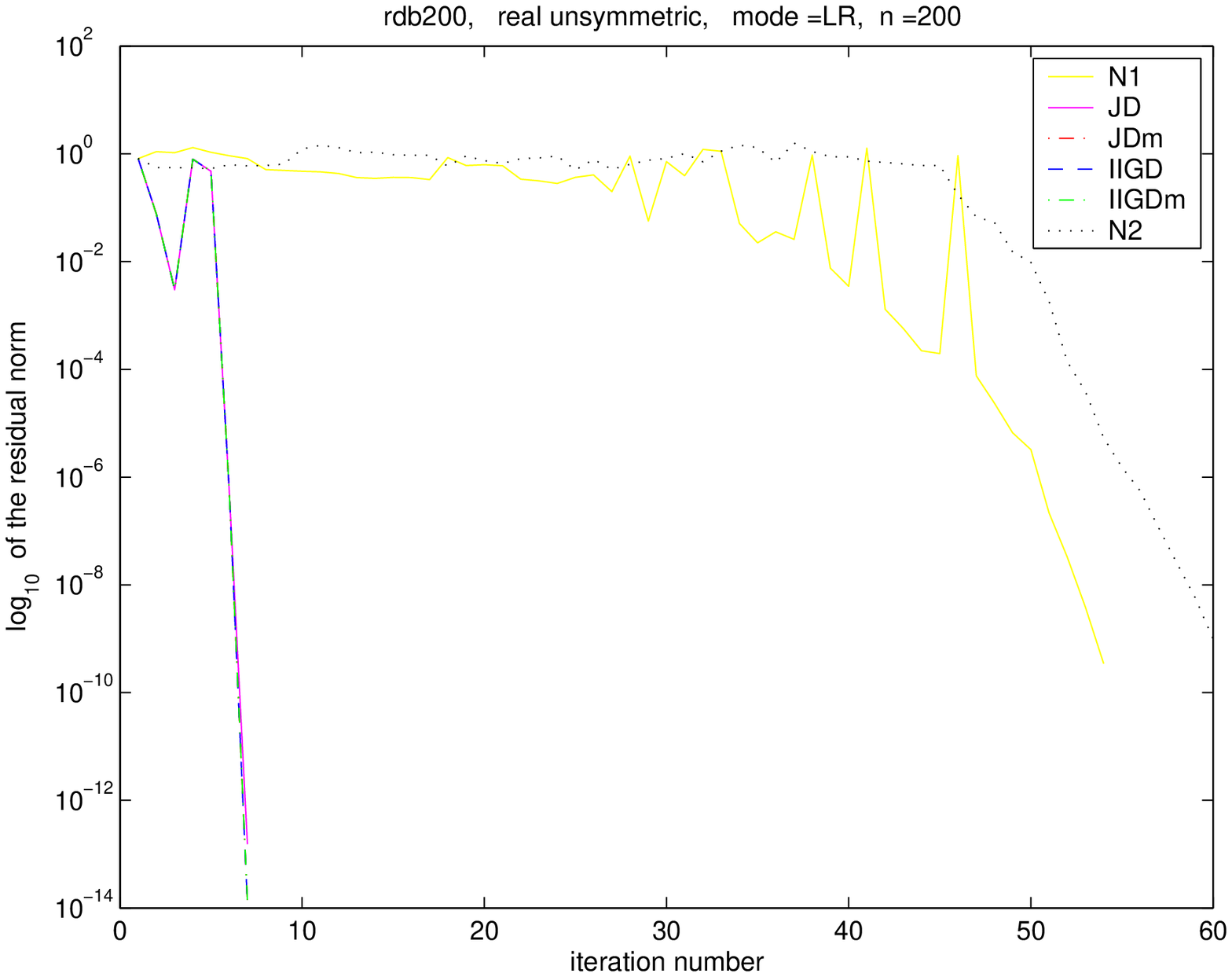}}
\scalebox{.397}{\includegraphics{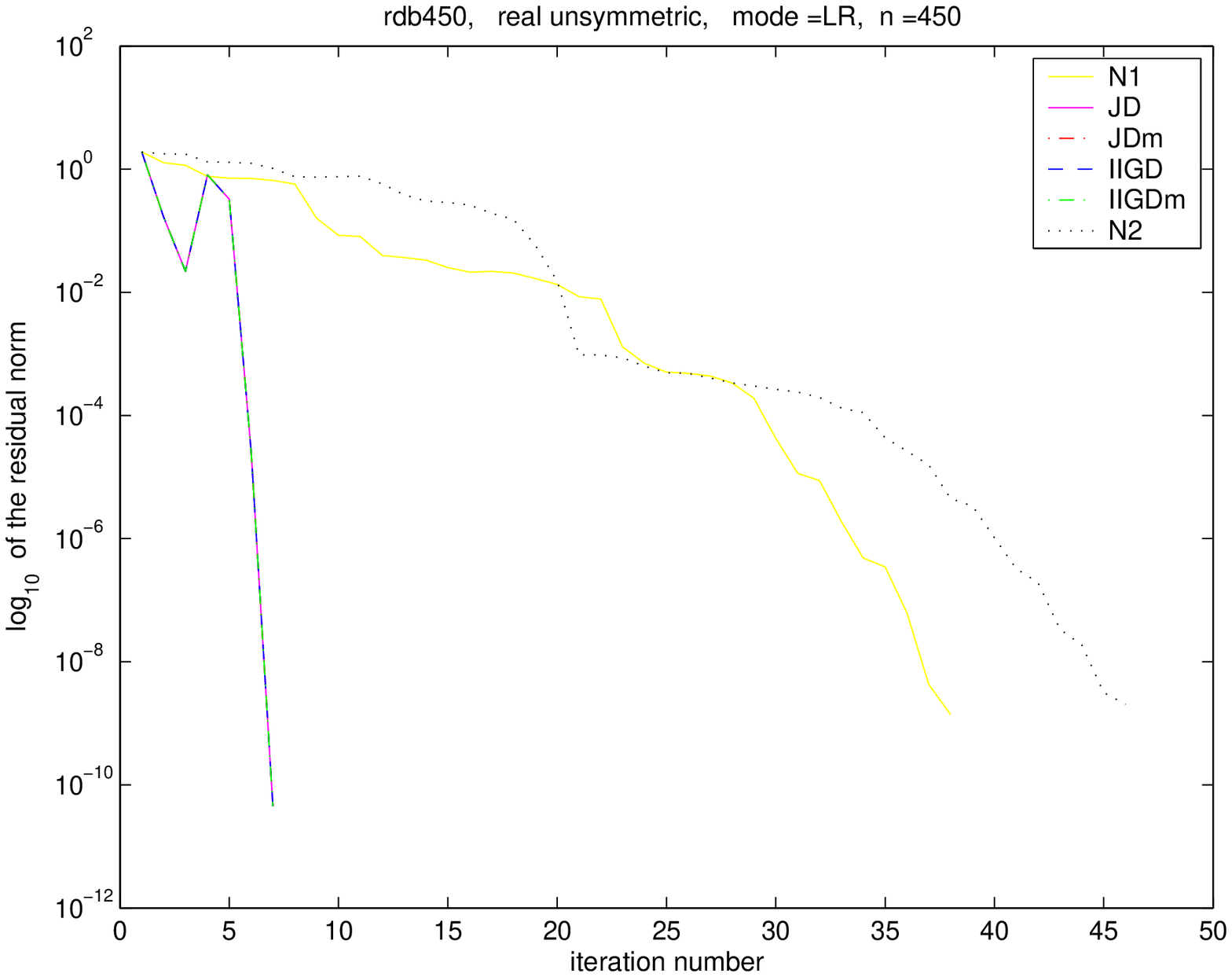}}
\end{center}

\caption{Comparison of performance for nonsymmetric models: 
rdb968, rw496, rdb200, rdb450. }
\label{nonsym2}

\end{figure}

\begin{figure}[htp]
\begin{center}
\scalebox{.397}{\includegraphics{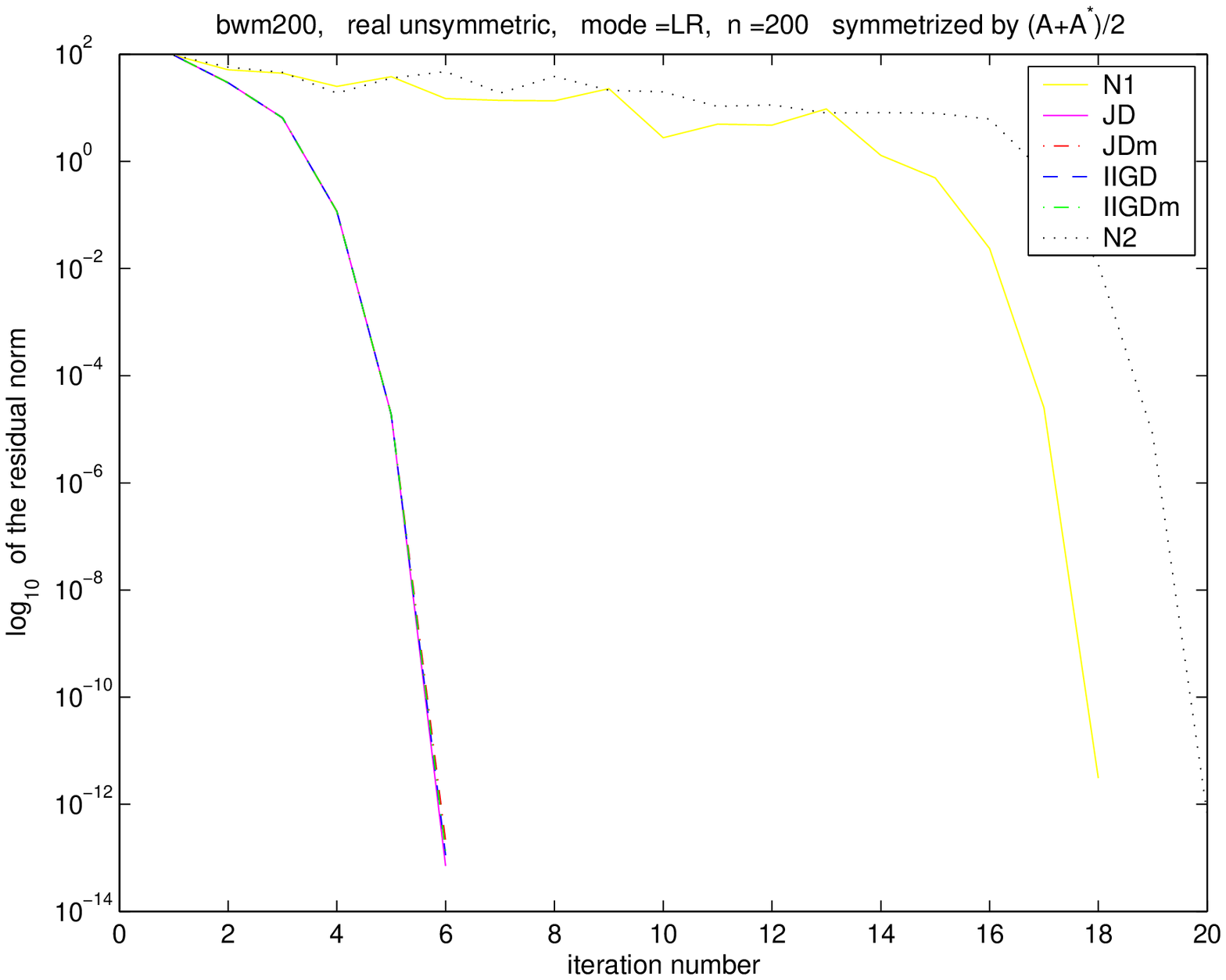}}
\scalebox{.397}{\includegraphics{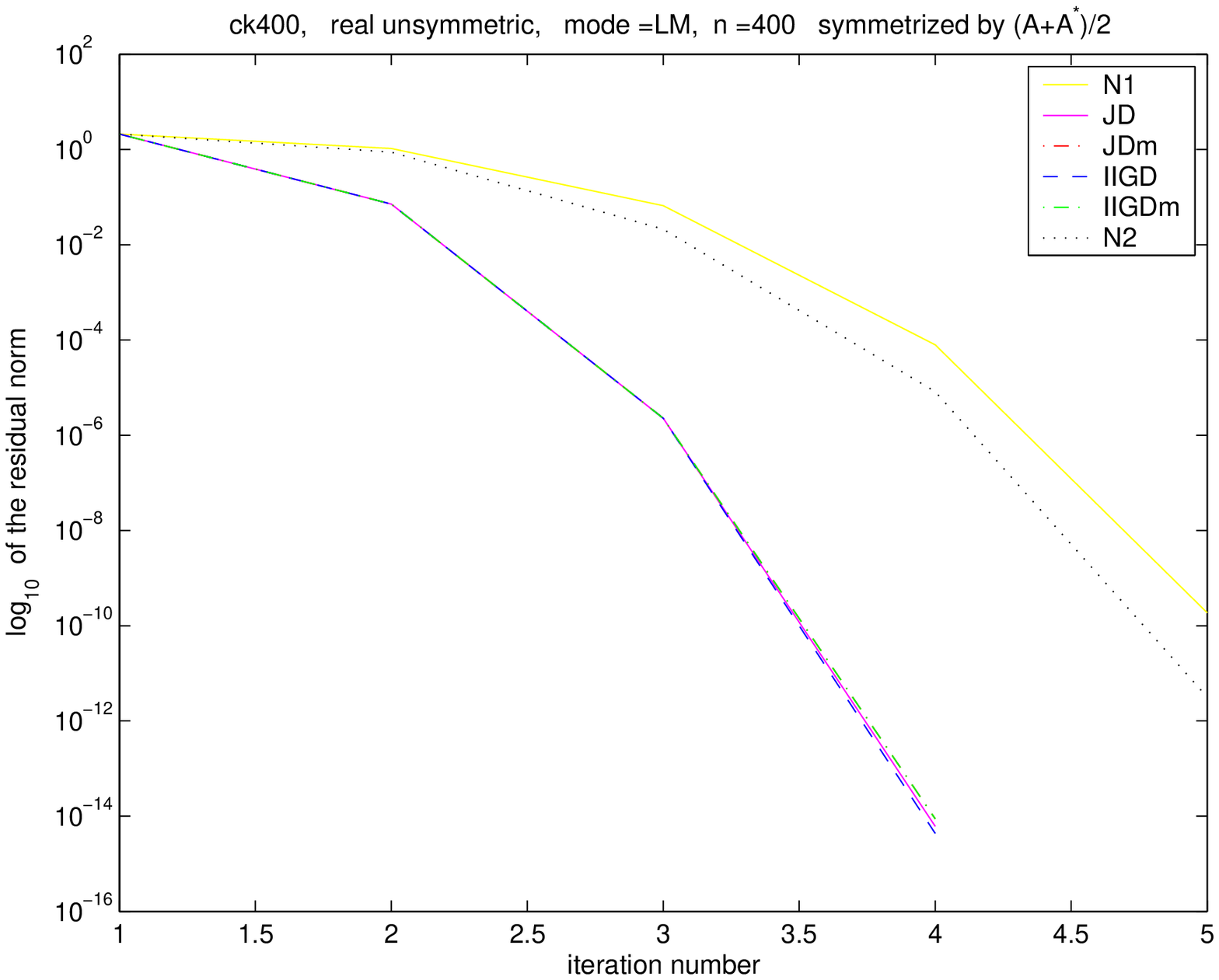}}

\vspace*{.2cm}

\scalebox{.397}{\includegraphics{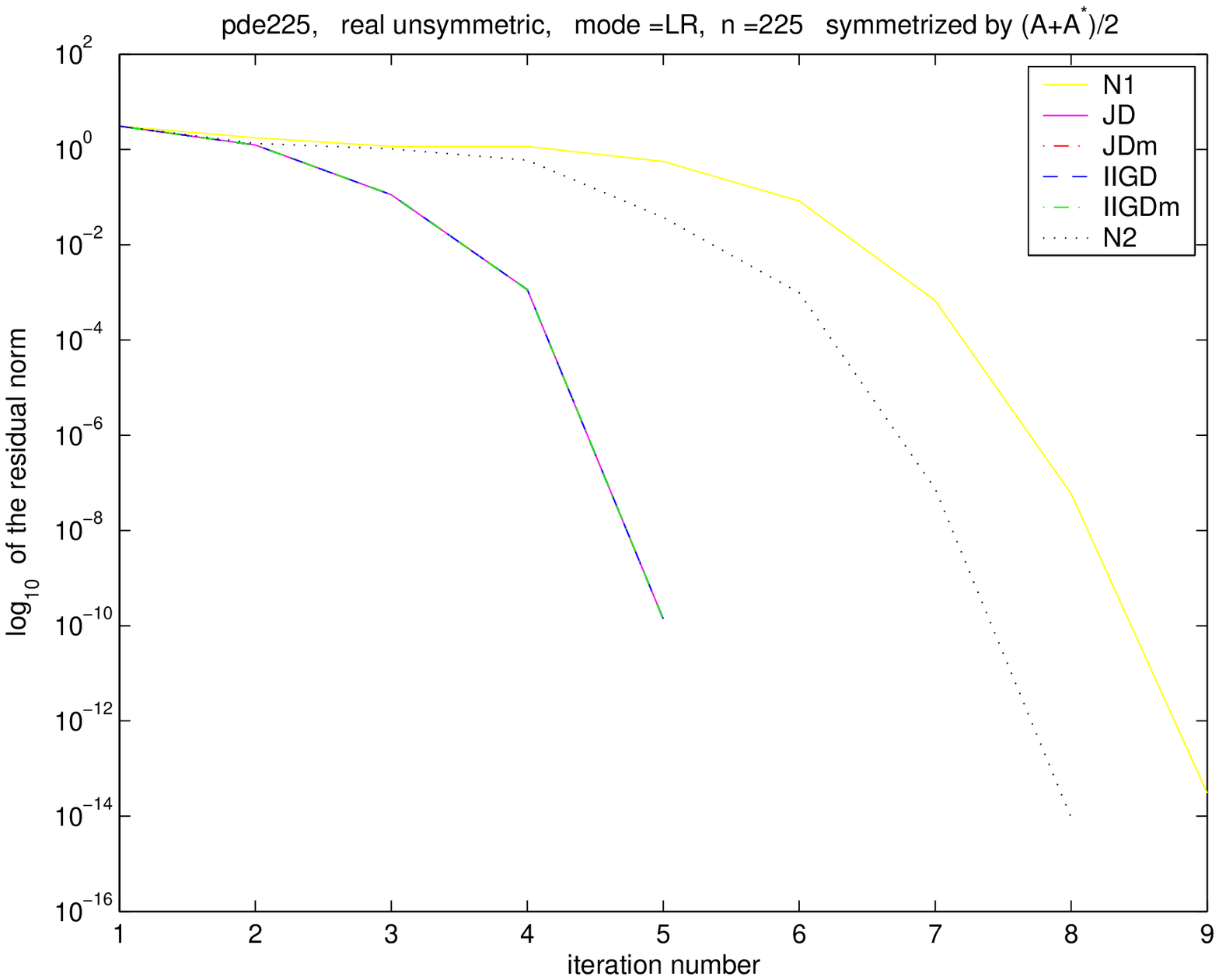}}
\scalebox{.397}{\includegraphics{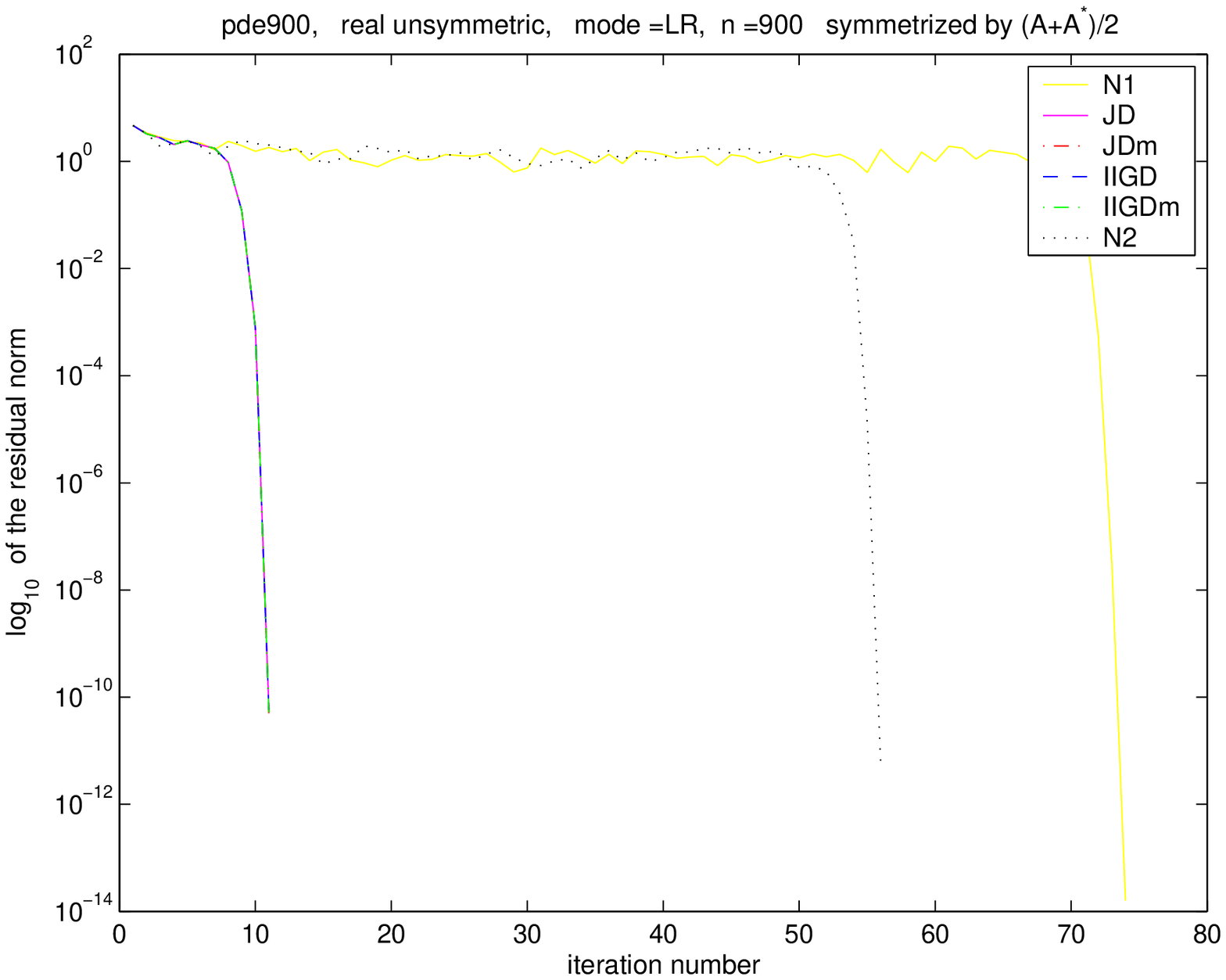}}

\vspace*{.2cm}

\scalebox{.397}{\includegraphics{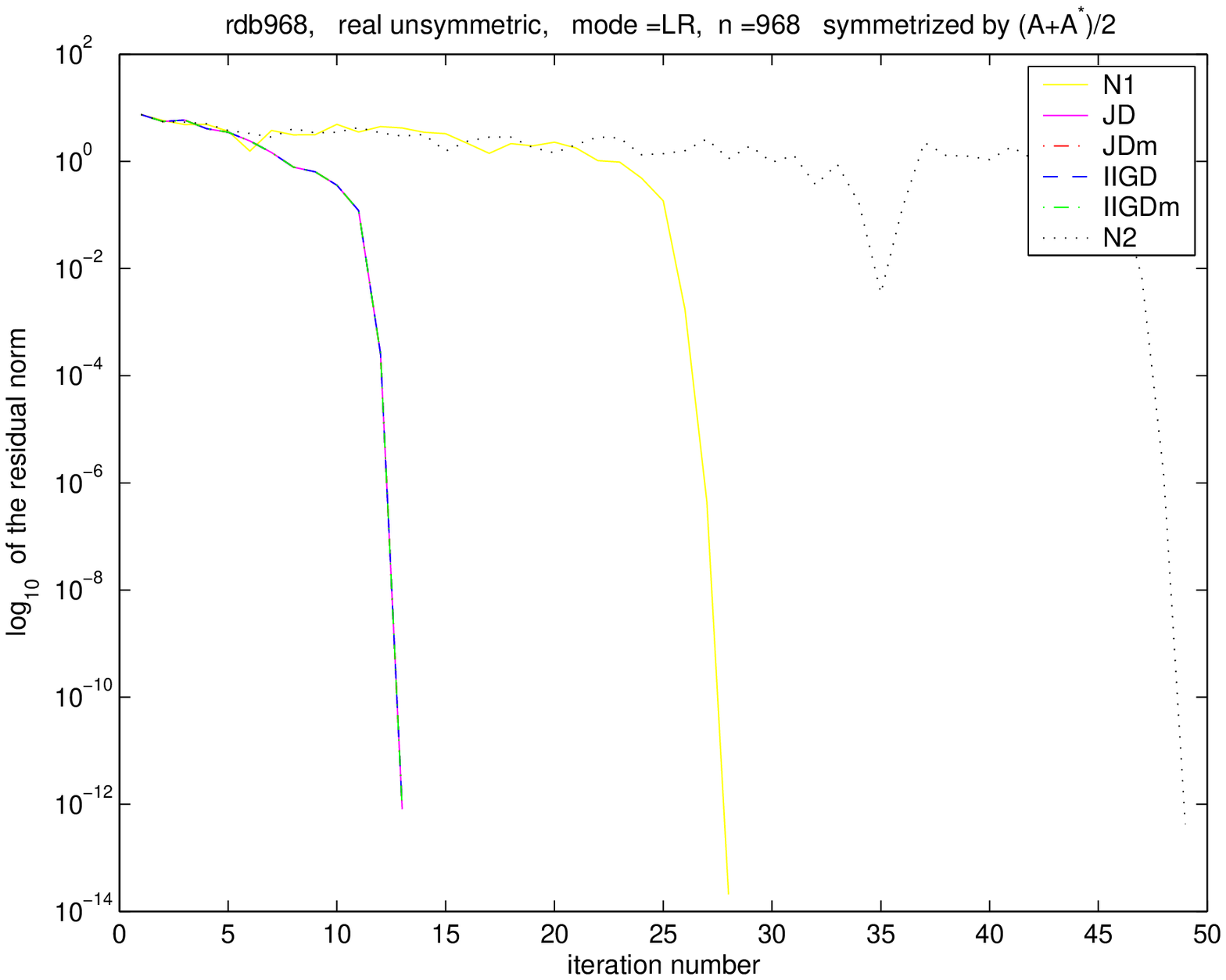}}
\scalebox{.397}{\includegraphics{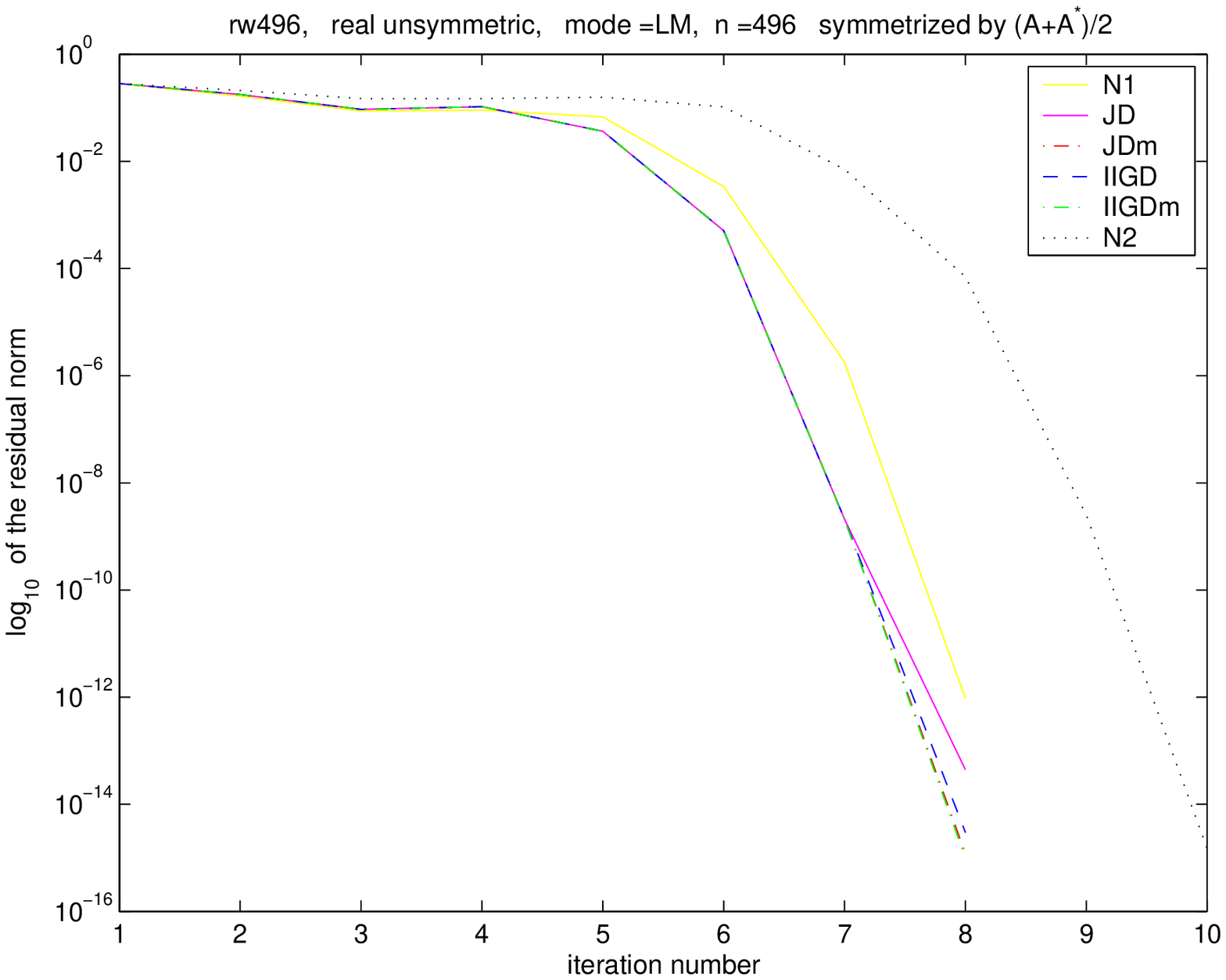}}

\end{center}
\caption{Comparison of performance for symmetric models: bwm200s, ck400s, pde225s, 
pde900s,  rdb968s, rw496s.}
\label{sym1}
\end{figure}

\begin{figure}[ht]
\begin{center}

\scalebox{.397}{\includegraphics{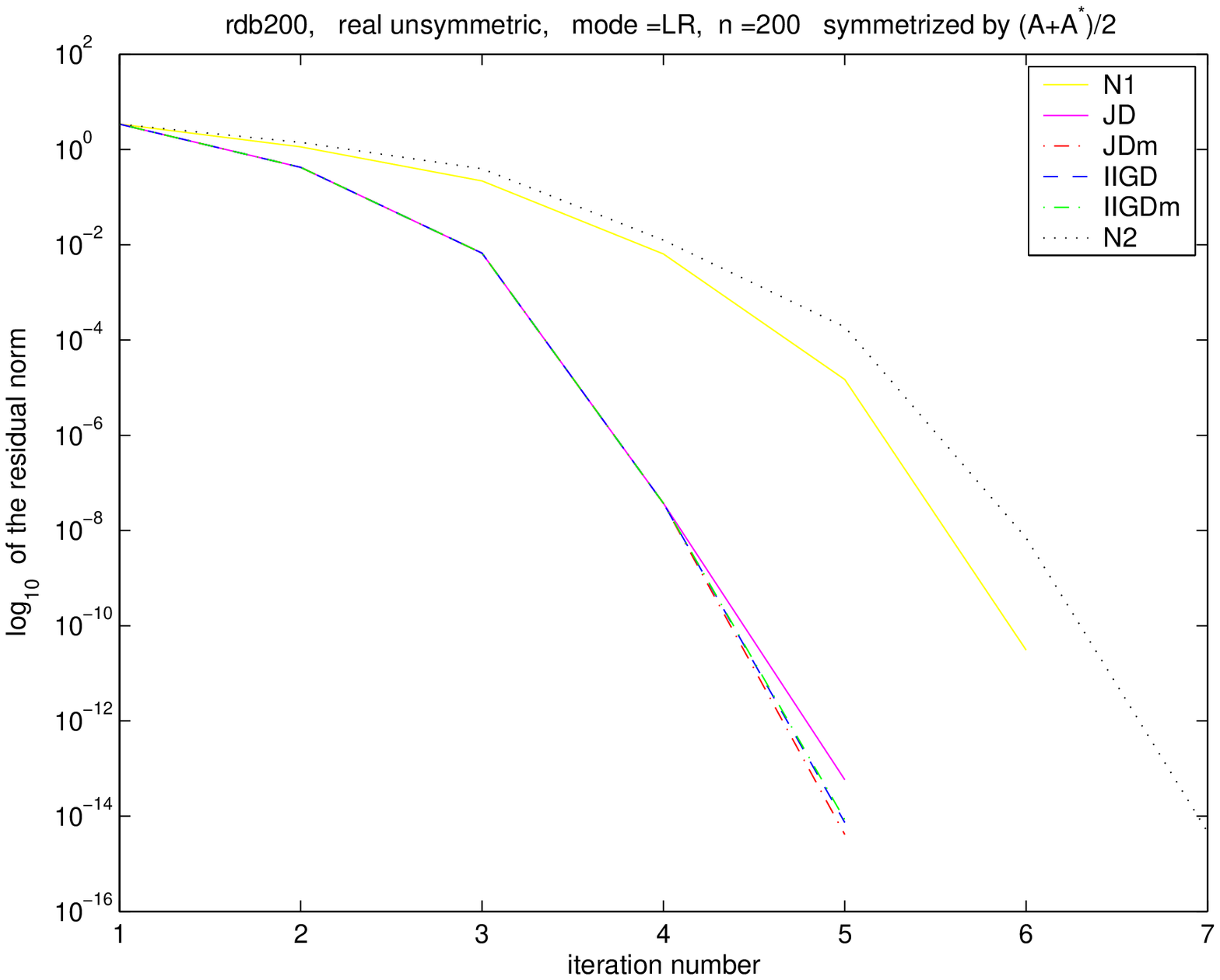}}
\scalebox{.397}{\includegraphics{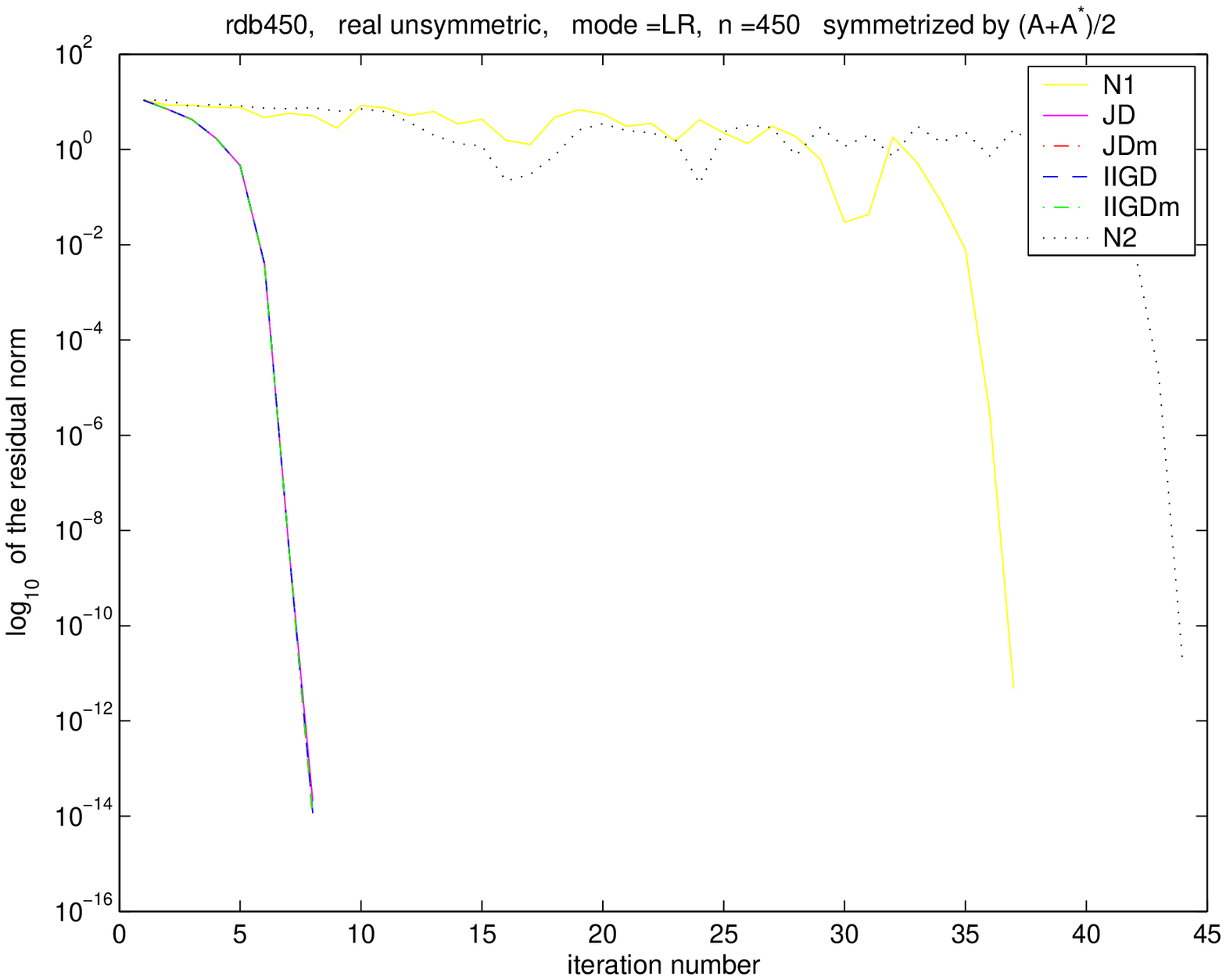}}

\vspace*{.2cm}
\scalebox{.396}{\includegraphics{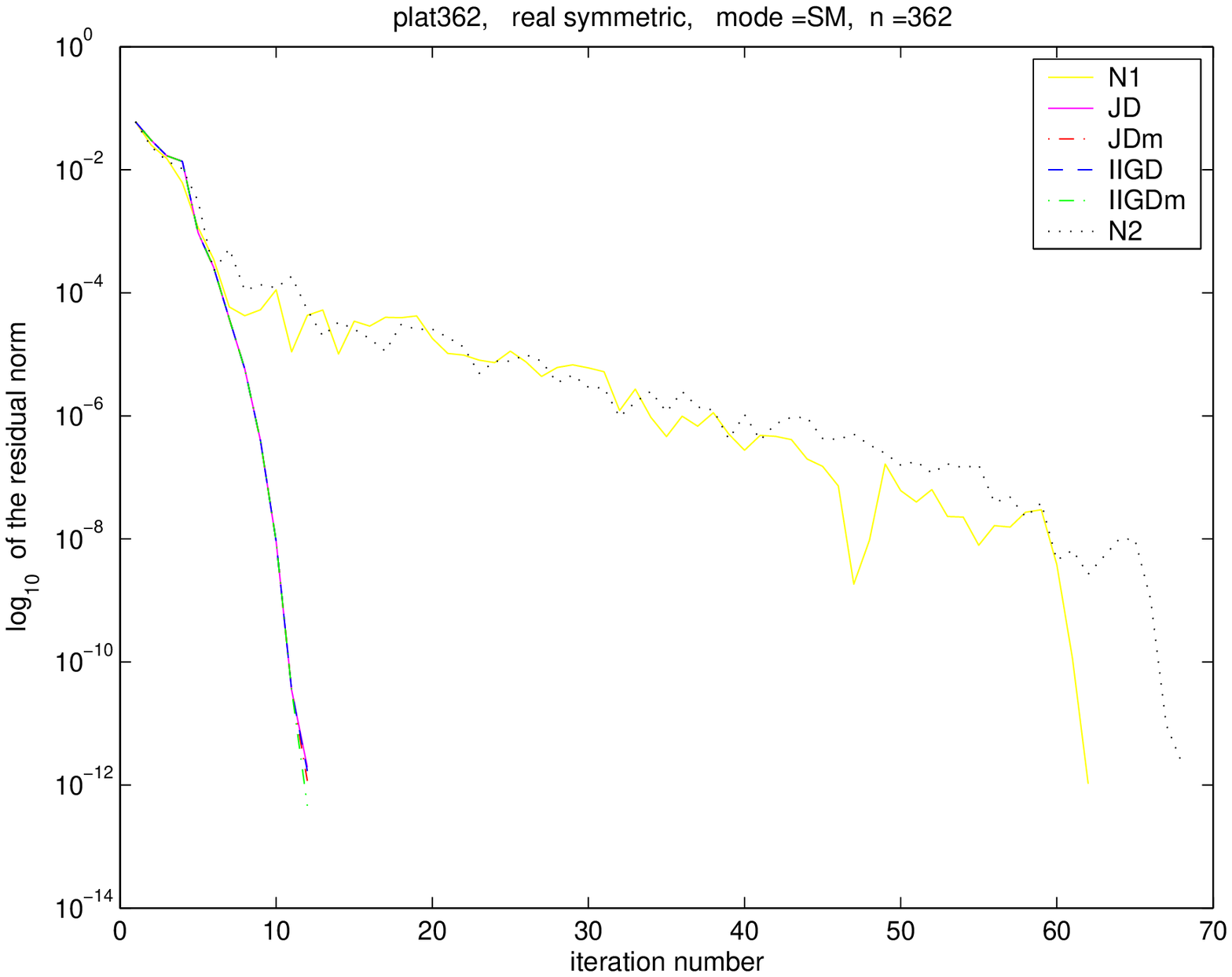}}
\scalebox{.396}{\includegraphics{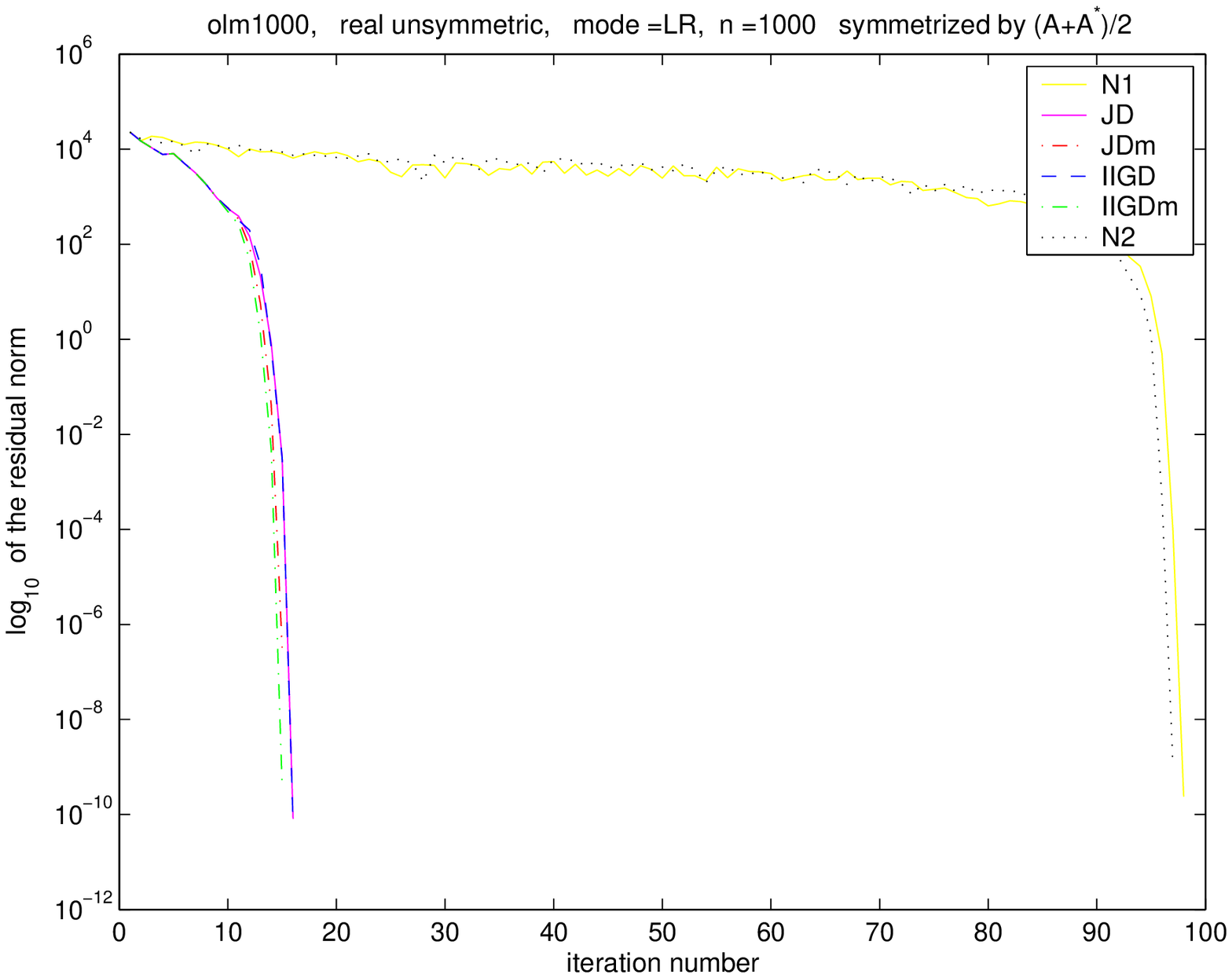}}
\end{center}

\caption{Comparison of performance for symmetric models:  rdb200s, rdb450s,
 plat362, olm1000s. }
\label{sym2}

\end{figure}

\section{Concluding remarks}

In this paper we highlighted the value of subspace methods in eigenvalue computation.
We also emphasized the link between eigenvalue computation and optimization through the
Newton method.
We derived a class of Newton updates for symmetric eigenvalue problems. 
From the study of the Newton updates and their similarity to the Jacobi-Davidson method,
we discovered the generalized correction equation \rf{tmp2}.  From this generalized 
formula we arrived at simplified forms of the Jacobi-Davidson method 
and the IIGD method.  The 
modification \rf{rqieq} actually halves the computational cost for the inner
iteration of IIGD.  We showed by extensive numerical results that the simplified
versions perform as efficiently as the Jacobi-Davidson and IIGD methods.  Preconditioned 
forms based on \rf{jdmeq} and \rf{rqieq} are expected to be efficient and practical
for large scale eigenvalue problems from real-world applications; this will be the subject
of ongoing research. \\

{\bf Acknowledgments:}
{The author thanks Prof. Dan Sorensen for introducing him to eigenvalue
computations. He thanks 
R. Shepard and M. Minkoff for some interesting discussions. The valuable NSF
travel support from the Institute of Pure and Applied Mathematics (IPAM) to 
attend the Nanoscale science and Engineering workshops at UCLA is also acknowledged.}


\begin{thebibliography}{ZZ}

\bibitem{arnoldi} W. E. Arnoldi, 
   {\it The Principle of Minimized Iterations in the Solution
    of the Matrix Eigenvalue Problem},
    Quarterly Journal of Applied Mathematics, 9 (1951) pp. 17-29.

\bibitem{temp00} Z. Bai, J. Demmel, J. Dongarra,
                  A. Ruhe, H. van der Vorst (Ed.),  
   {\it Templates for the Solution of Algebraic
                  Eigenvalue Problems: A Practical Guide}, 
    SIAM publications, 2000.

\bibitem{chatel93} F. Chaitin-Chatelin, 
   {\it Eigenvalues of Matrices},  John Wiley \& Son, 1993.

\bibitem{cps94} M. Crouzeix, B. Philippe, M. Sadkane,
   {\it The Davidson Method,}  SIAM Journal on Scientific Computing, 15 (1994) pp. 62-76.

\bibitem{cul_wil02} J. K. Cullum,  R. A. Willoughby,
   {\it Lanczos Algorithms for Large Symmetric Eigenvalue Computations: Vol. I: Theory},
   SIAM, Classics in Applied Mathematics 41, 2002.

\bibitem{dgks} J. Daniel, W. B. Gragg, L. Kaufman, and G. W. Stewart,
   {\it Reorthogonalization and Stable Algorithms for
    Updating the {Gram}-{Schmidt} {QR} Factorization},
    Mathematics of Computation, 30 (1976) pp. 772-795.

\bibitem{davids75} E. R. Davidson,
   {\it The Iterative Calculation of A Few of the Lowest Eigenvalues
    and Corresponding Eigenvectors of Large Real Symmetric Matrices,}
    Journal of Computational Physics, 17 (1975) pp. 87-94.

\bibitem{davids93} E. R. Davidson,
   {\it Monster Matrices: Their Eigenvalues and Eigenvectors},
   Computers in Physics, 7 (1993) pp. 519-522.

\bibitem{dax03} A. Dax,
   {\it The Orthogonal Rayleigh Quotient Iteration Method},
   Linear Algebra and Its Applications, 358 (2003) pp. 23-43.

\bibitem{den_sch83} J. E. Dennis and R. Schnabel,
   {\it Numerical Methods for Unconstrained Optimization and
     Nonlinear Equations}, Prentice-Hall, 1983; updated reprint,
   SIAM classics in Applied Mathematics, 16, 1996.

\bibitem{ddsv98} J. Dongarra, I. Duff, D. C. Sorensen, H. A. van der Vorst,
   {\it Numerical Linear Algebra for High-Performance Computers},
   SIAM publications, 1998.


\bibitem{fsv98} D. R. Fokkema, G. L. G. Sleijpen, H. A. van der Vorst,
   {\it Jacobi-Davidson Style QR and QZ Algorithms for the Reduction of Matrix Pencils,}
   SIAM Journal on Scientific Computing, 20 (1998) pp. 94-125.


\bibitem{gv96} G. H. Golub and C. F. Van Loan,
   {\it Matrix Computations},
   Johns Hopkins University Press, 3rd ed., 1996.

\bibitem{gol_ye02} G. H. Golub, Q. Ye,
   {\it An Inverse Free Preconditioned Krylov Subspace Method for Symmetric 
    Generalized Eigenvalue Problems},
    SIAM Journal on Scientific Computing, 24 (2002) pp. 312-334.

\bibitem{greenb97} A. Greenbaum,
   {\it Iterative Methods for Solving Linear Systems},
   SIAM,  Frontiers in Applied Mathematics, v17, 1997.

\bibitem{hoc_sle03} M. E. Hochstenbach, G. L. G. Sleijpen,
   {\it Two-sided and Alternating Jacobi-Davidson,}
   Linear Algebra and Its Applications, 358 (2003) pp. 145-172.


\bibitem{hor_joh91} R. A. Horn and C. R. Johnson,
   {\it Topics in Matrix Analysis},
    Cambridge University Press, 1991.


\bibitem{knya01} A. V. Knyazev,
   {\it Toward the Optimal Preconditioned Eigensolver: Locally Optimal
   Block Preconditioned Conjugate Gradient Method,}
   SIAM Journal on Scientific Computing, 23 (2001) pp. 517-541. 
\bibitem{kny_ney03} A. V. Knyazev, Klaus Neymeyr,
  {\it Efficient Solution of Symmetric Eigenvalue Problems Using Multigrid
   Preconditioners in the Locally Optimal Block Conjugate Gradient Method,}
   Electronic Transactions on Numerical Analysis, 15 (2003), pp. 38-55. 

\bibitem{lanczos} C. Lanczos,
   {\it An Iteration Method for the Solution of the Eigenvalue Problem
    of Linear Differential and Integral Operators},
    J. Res. Nat. Bur. Standards, 45 (1950) pp. 255-282.

\bibitem{leh_mee98} R. B. Lehoucq, K. Meerbergen,
   {\it Using Generalized Cayley Transformations within an Inexact Rational
    Krylov Sequence Method,}
   SIAM Journal on Matrix Analysis and Applications, 20 (1998) pp. 131-148.

\bibitem{mor_sco86} R. B. Morgan, D. S. Scott,
   {\it Generalization of Davidson's Method for Computing Eigenvalues
   of Sparse Symmetric Matrices}, SIAM Journal on Statistical and Scientific Computing, 
   7 (1986) pp. 817-825.  


\bibitem{notay03} Y. Notay,
  {\it Convergence Analysis of Inexact Rayleigh Quotient Iteration,}
  SIAM Journal on Matrix Analysis and Applications, 24 (2003) pp. 627-644.

\bibitem{ojs90} J. Olsen, P. J{\o}rgensen, J. Simons,
   {\it Passing the One-billion Limit in Full
   Configuration-Interaction (FCI) Calculations,} Chemical Physics Letters,
   169 (1990) pp. 463-472.

\bibitem{ostrow} A. M. Ostrowski, 
   {\it On the Convergence of the Rayleigh Quotient Iteration for the Computation
    of Characteristic Roots and Vectors. I--VI,} Arch. Rational Meth. Anal., v. 1--4,
    1958/59, pp. 233-241,423-428,325-340,341-347,472-481,153-165.
    (cited in \cite{parlet74}). 


\bibitem{parlet74} B. N. Parlett,
   {\it The Rayleigh Quotient Iteration and Some Generalizations for Nonnormal Matrices,}
   Mathematics of Computation, 28 (1974) pp. 679-693.

\bibitem{parlet80} B. N. Parlett,
   {\it The Symmetric Eigenvalue Problem}, Prentice-Hall, 1980; updated reprint, 
    SIAM, Classics in Applied Mathematics 20, 1998.

\bibitem{pet_wil79} G. Peters, J. H. Wilkinson,
   {\it Inverse Iteration, Ill-conditioned Equations and Newton's Method},
   SIAM Review, 21 (1979) pp. 339-360.

\bibitem{ruhe84} A. Ruhe,
   {\it Rational Krylov Sequence Methods for Eigenvalue Computations},
   Linear Algebra and Its Applications, 58 (1984) pp. 391-405.

\bibitem{saad_eig92} Y. Saad,
   {\it Numerical Methods for Large Eigenvalue Problems,}
   John Wiley, 1992.  \\
    \url{http://www-users.cs.umn.edu/~saad/books.html}.




\bibitem{saad86} Y. Saad and M. Schultz,
   {\it {GMRES}: A Generalized Minimal Residual Algorithm for
                 Solving Nonsymmetric Linear Systems},
   SIAM Journal on Statistical and Scientific Computing, 7 (1986) pp. 856-869.


\bibitem{jd96} G. L. G. Sleijpen and H. A. van der Vorst, 
   {\it A Jacobi-Davidson Iteration Method for Linear Eigenvalue Problems}, 
   SIAM Journal on Matrix Analysis and Applications,  17 (1996) pp. 401-425.  
   Reprinted in SIAM Review, 42 (2000), pp. 267-293.

\bibitem{gjd96}G. L. G. Sleijpen, A. G. L. Booten, D. R. Fokkema, H. A. van der Vorst,
   {\it  Jacobi-Davidson Type Method for Generalized Eigenproblems and
    Polynomial Eigenproblems,} BIT, 36 (1996) pp. 595-633.


\bibitem{jd_in} G. L. G. Sleijpen and H. A. van der Vorst,
   {\it The Jacobi-Davidson Method for Eigenvalue Problems and Its Relation 
   to Accelerated Inexact Newton Schemes,}
   Proceeding of the second IMACS International Symposium on 
   Iterative Methods in Linear Algebra, June, 1995.


\bibitem{sorens92} D. C. Sorensen,
   {\it Implicit Application of Polynomial Filters
   in a k-step Arnoldi Method},
   SIAM Journal on Matrix Analysis and Applications,  13 (1992) pp. 357-385.

\bibitem{sorens02} D. C. Sorensen,
   {\it Numerical Methods for Large Eigenvalues Problems},
   Acta Numerica,  (2002)  pp. 519-584.

\bibitem{sor_yan98} D. C. Sorensen, C. Yang,
   {\it A Truncated RQ Iteration for Large Scale Eigenvalue Calculations,}
   SIAM Journal on Matrix Analysis and Applications, 19 (1998) pp. 1045-1073.

\bibitem{ssf95} A. Stathopoulos, Y. Saad, C. F. Fisher,
   {\it Robust Preconditioning of Large, Sparse, Symmetric Eigenvalue Problems,}
   Journal of Computational and Applied Mathematics, 64 (1995) pp. 197-215. 

\bibitem{ssw98} A. Stathopoulos, Y. Saad, K. Wu,
   {\it Dynamic Thick Restarting of the Davidson, and the Implicit
   Restarted Arnoldi Methods,}
   SIAM Journal on Scientific Computing,  19 (1998) pp. 227-245.

\bibitem{stewar00} G. W. Stewart,
   {\it Matrix Algorithms,  Volume II: Eigensystems,}
   SIAM publications, 2001.
\bibitem{stewar01}  G. W. Stewart,
  {\it A Krylov--Schur Algorithm for Large Eigenproblems},
    SIAM Journal on Matrix Analysis and Applications, 23 (2001) pp. 601-614.



\bibitem{wilkin65} J. H. Wilkinson, 
   {\it The Algebraic Eigenvalue Problem,}
    Oxford University Press, 1965.


\bibitem{wss98} K. Wu, Y. Saad, A. Stathopoulos,
   {\it Inexact Newton Preconditioning Techniques for Large
              Symmetric Eigenvalue Problems},
    Electronic Transactions on Numerical Analysis, 7 (1998) pp. 202-214.

\bibitem{wu_sim98} K. Wu, H. Simon,
   {\it Thick-restart Lanczos Method for Large Symmetric Eigenvalue Problems,}
    SIAM Journal on Matrix Analysis and Applications, 22 (2000) pp. 602-616.

\bibitem{yang98} C. Yang,
  {\it Convergence Analysis of an Inexact Truncated RQ Iterations,}
   Electronic Transactions on Numerical Analysis, 7 (1998) pp. 40-55.

\bibitem{bicgstab} H. A. van der Vorst,
   {\it     Bi-CGSTAB: A Fast and Smoothly Converging Variant of Bi-CG
    for the Solution of Nonsymmetric Linear Systems,}
   SIAM Journal on Statistical and Scientific Computing, 13 (1992) pp. 631-644.



\end{thebibliography}
\end{document}